\documentclass[a4paper, 11pt ]{amsart}
\usepackage{mathrsfs,amssymb}
\usepackage{multicol}\multicolsep=0pt

\usepackage[sort]{cite}

\def\crulefill{\leavevmode\leaders\hrule height 1pt\hfill\kern 0pt}
\long\def\QUERY#1{%
\leavevmode\newline%
\noindent$\star\star\star$\thinspace\textsf{Comment/Query}\crulefill\newline%
   \space #1\newline\hbox to 120mm{\crulefill}$\star\star\star$\newline}

\numberwithin{equation}{section} \theoremstyle{definition}
\newtheorem{Defn}[equation]{Definition}
\newtheorem{Example}[equation]{Example}
\newtheorem{Remark}[equation]{Remark}
\theoremstyle{plain}
\newtheorem{Prop}[equation]{Proposition}
\newtheorem{Theorem}[equation]{Theorem}
\newtheorem{Assumption}[equation]{Assumption}
\newtheorem{Lemma}[equation]{Lemma}
\newtheorem{Cor}[equation]{Corollary}

\def\Case#1{\medskip\noindent\textbf{Case #1}:\leavevmode\newline}

\makeatletter
\def\enumerate{\begingroup\ifnum\@enumdepth>3\@toodeep\else
      \advance\@enumdepth\@ne
      \edef\@enumctr{enum\romannumeral\the\@enumdepth}%
      \topsep\z@\parskip\z@
      \list{\csname label\@enumctr\endcsname}
        {\@nmbrlisttrue\let\@listctr\@enumctr
         \parsep\z@\itemsep\z@\topsep\z@
         \setcounter{\@enumctr}{0}
         \def\makelabel##1{\hss\llap{\rm ##1}}
       }\fi}

\makeatother

\let\epsilon=\varepsilon
\def\({\big(}
\def\){\big)}

\def\R{\mathbb R}

\def\0{\underline{0}}

\def\Sym{\mathfrak S}

\DeclareMathOperator{\rank}{rank}

\def\simk{\overset{k}\sim}

\def\simnn{\overset{n-1}\sim}

\def\half{\frac12}

\let\gdom\rhd
\let\gedom\unrhd

\def\m{\mathfrak m}
\def\floor#1{\lfloor\tfrac#1\rfloor}
\def\UPD{\mathscr{T}^{ud}}
\def\Std{\mathscr{T}^{std}}

\def\m{\mathfrak m}
\def\n{\mathfrak n}
\def\s{\mathfrak s}
\def\ts{\tilde\s}
\def\t{\mathfrak t}
\def\u{\mathfrak u}
\def\v{\mathfrak v}

\def\ss{\mathbf s}
\def\ts{\mathbf t}
\def\us{\mathbf u}

\newcommand{\cba}[1]{\mathscr {B}_{#1}({\delta})}

\newcommand{\ccba}[1]{\mathscr {B}_{#1}({0})}

{\catcode`\|=\active
  \gdef\set#1{\mathinner{\lbrace\,{\mathcode`\|"8000%
                                   \let|\midvert #1}\,\rbrace}}
  \gdef\seT#1{\mathinner{\Big\lbrace\,{\mathcode`\|"8000%
                                   \let|\midverT #1}\,\Big\rbrace}}
}
\def\midvert{\egroup\mid\bgroup}
\def\midverT{\egroup\,\Big|\,\bgroup}

\def\Set[#1]#2|#3|{\Big\{\ #2\ \Big| \
           \vcenter{\hsize #1mm\centering #3}\Big\}}

\def\map#1#2{\,{:}\,#1\!\longrightarrow\!#2}

\title{Discriminants of Brauer algebras}
\author{Hebing Rui and Mei Si}
\address{H.R. Department of Mathematics,  East China Normal
University, Shanghai, 200062, China}
\email{hbrui@math.ecnu.edu.cn}
\address{M. S. Department of Mathematics,  East China Normal
University, Shanghai, 200062, China}
\email{52050601011@student.ecnu.edu.cn}
\thanks{The first author is supported in part by NSFC 
and NCET}

\begin{document}
\baselineskip15pt
\begin{abstract} In this paper, we compute
Gram determinants associated to all  cell modules of  Brauer
algebras $\cba{n}$.  Theoretically, we know when a cell  module of
$\cba{n}$ is equal to its simple head. This gives a solution of this
long standing problem.
\end{abstract}

\sloppy \maketitle

\centerline{\textit{On the occasion of Professor Gus Lehrer's
                      60$^{\text{th}}$ birthday}}

\section{Introduction}
In~\cite{BrauerAlg}, Richard Brauer  introduced a  class of finite
dimensional algebras $\cba{n}$  in order to study the $n$-th tensor
power of the defining representations of orthogonal groups and
symplectic groups. Such algebras, which are known as  Brauer
algebras or Brauer centraliser algebras,  have been studied by many
authors; see for example\cite{BM,Brown:brauer,DoranHanlonWales,
Enyang1, Enyang,
FiGr,HalRam:basic,HW2,HW1,J:Brauer,LeducRam,Nazarov:brauer,Ram:chBrauer,Rui:ssbrauer,
RS:ssbrauer,Terada:brauer,Wenzl:ssbrauer,Weyl}, etc.

In \cite{GL}, Graham and Lehrer  introduced the notion of cellular
algebra which is defined over a poset $\Lambda$. Such an algebra has
a nice basis, called a cellular basis. For each $\lambda\in
\Lambda$, there is a module  $\Delta(\lambda)$ called a cell module.
Graham and Lehrer  showed  that there is a symmetric, associative
bilinear form $\phi_{\lambda}$ defined on $\Delta(\lambda)$. It has
been proved in \cite[3.8]{GL} that a cellular algebra is (split)
semisimple if and only if $\phi_{\lambda}$ is non-degenerate for any
$\lambda\in \Lambda$.

Graham and Lehrer  proved that $\cba{n}$ over a commutative ring is
a cellular algebra over the poset $\Lambda$ which consists of all
pairs $(f, \lambda)$, with $0\le f\le \lfloor n/2\rfloor$ and
$\lambda$ being a partition of $n-2f$. Here $\lfloor n/2\rfloor$ is
the maximal integer with $\lfloor n/2\rfloor\le n/2$. Therefore,  $\cba{n}$
 is (split) semisimple over a field if and only if $\phi_{f,
\lambda}$ are non-degenerate for all $(f, \lambda)\in \Lambda$.

 In \cite{Rui:ssbrauer}, the first author  proved that
 $\mathscr B_n(\delta)$ with  $\delta\neq 0$ is semisimple if and only if (a) 
 the group algebra of the symmetric group
 $\mathfrak S_n$ is semisimple, (b) 
$\phi_{1, \lambda}$ are non-degenerate for all partitions $\lambda$ of
$k-2$ and $2\le k\le n$.  By  \cite[3.3-3.4]{DoranHanlonWales}, 
he gave an algorithm to determine whether $\phi_{1, \lambda}$ is 
non-degenerate or not. When $\delta=0$, the semisimplicity of $\mathscr B_n(\delta)$ 
can be determined by direct computation \cite{Rui:ssbrauer}. 
 In \cite{RS:ssbrauer}, we used such an algorithm to verify Enyang's  
 conjecture in \cite{Enyang}, which is about the
semisimplicity of Brauer algebras. This gives a complete solution of
the problem of semisimplicity of $\mathscr B_n(\delta)$ over an
arbitrary field.  In other words we have   found all $\delta$ such
that no  Gram determinant associated to a cell module is equal to zero. 
 
A further question is when the Gram determinant for a fixed cell
module is not equal to zero. This is equivalent to finding  a
necessary and sufficient condition for a cell module of $\cba{n}$
being equal to its simple head.

The main purpose of this  paper is to solve the above  problem. In
\cite{Enyang}, Enyang  constructed the Jucys-Murphy basis of each
cell module for $\cba{n}$.  In this paper, we will use his result to
construct the Jucys-Murphy basis  for $\cba{n}$ over a commutative
ring.
 Since the Jucys-Murphy elements of $\cba{n}$, which were  defined by Nazarov
in \cite{Nazarov: brauer},  act on the Jucys-Murphy basis of
$\cba{n}$ via upper triangular matrices, we can and do use some
arguments in \cite{JM} or \cite{M:gendeg} to construct an orthogonal
basis of $\cba{n}$ under the assumption~\ref{assumption}. We show
that  the Gram determinants associated to each cell module, which
are defined via the Jucys-Murphy basis and this orthogonal basis, are same.
Using classical branching rule for $\cba{n}$, we will compute each
diagonal entry of the Gram matrix defined via this  orthogonal basis. This
gives a  recursive formula for each Gram determinant. We remark that
such recursive formulae   shed some light on possible composition
factors of cell modules when they are not simple. Also, we will use
such  recursive formulae to describe (up to a sign) the actions of
generators of $\cba{n}$ on each orthogonal basis element. As an
application, we construct explicitly
 primitive idempotents and central primitive idempotents  of
$\cba{n}$. This gives  the Wedderburn-Artin decomposition of Brauer
algebras $\cba{n}$ under the assumption~\ref{assumption}.

Finally, we remark that the method in this paper will be used to
deal with Birman-Murakami-Wenzl algebras~\cite{BM} and cyclotomic
Nazarov-Wenzl algebras \cite{NW}.  Details will appear elsewhere.

\textbf{Acknowledgement:} The first author thanks Prof. A.~Mathas for
inviting him to the University of Sydney in January 2005,  and
Prof. G. Lehrer for suggesting to investigate discriminants for
Brauer algebras. We thank A. Mathas for his comments on the paper.
Finally, both of us thank the referee for their valuable comments.

\section{Brauer algebras}
Throughout, let $R$ be a commutative ring which contains the  identity $1$ and
$\delta$.  The Brauer algebra $\cba{n}$ is a unital associative
$R$-algebra generated by $s_i, e_i$, $1\le i\le n-1$, subject to the
following relations:
\begin{xalignat*}{3}
  s_i^2 &=1, & e_i^2 &=\delta e_i, & s_ie_i &=e_is_i=e_i,\\
s_is_j&=s_js_i,& s_ie_j&=e_js_i,& e_ie_j&=e_je_i,\\
  s_ks_{k+1}s_k&=s_{k+1}s_ks_{k+1},& e_ke_{k+1}e_k&=e_k,
                            & e_{k+1}e_ke_{k+1}&=e_{k+1},\\
  s_ke_{k+1}e_k&=s_{k+1}e_k, & e_{k+1}e_ks_{k+1}&=e_{k+1}s_k,
\end{xalignat*}
where $1\le i,j<n$, with $|i-j|>1$, and $1\le k<n-1$.

Let $Sym(T)$ be the symmetric group on the set $T$. If $T=\{1, 2,
\cdots, n\}$, we denote $Sym(T)$ by $\mathfrak S_n$. The group
algebra $R\mathfrak S_n$ can be considered as a subalgebra of
$\cba{n}$ if we identify $s_i$ with  the basic transposition $(i,
i+1)$.


We are going to construct the Jucys-Murphy basis of $\cba{n}$. We
start by recalling some combinatorics.

Recall that a partition of $n$ is a weakly decreasing sequence of
non--negative integers $\lambda=(\lambda_1,\lambda_2,\dots)$ such
that $|\lambda|:=\lambda_1+\lambda_2+\cdots=n$. In this situation,
we write $\lambda\vdash n$. Let $\Lambda^+(n)$ be the set of all
partitions of $n$. It is well known that $\Lambda^+(n)$ is a poset
with dominance order $\trianglelefteq$ as the partial order on it.
Given $\lambda, \mu\in \Lambda^+(n)$, $\lambda\trianglelefteq \mu$
if  $ \sum_{j=1}^i \lambda_j\le \sum_{j=1}^i \mu_j$ for all possible
$i$. Write $\lambda\vartriangleleft \mu$ if $\lambda\trianglelefteq
\mu$ and $\lambda\ne \mu$.

Suppose that  $\lambda$ and $\mu$ are two partitions. We say that
$\mu$ is obtained from $\lambda$ by adding a box if there exists an
$i$  such that $\mu_i=\lambda_i+1$ and $\mu_j=\lambda_j$ for all
$j\neq i$. In this situation we will also say that $\lambda$ is
obtained from $\mu$ by removing a box and we write
$\lambda\rightarrow\mu$ and $\mu\setminus\lambda=(i,\mu_i)$. We will
also say that the pair  $(i,\mu_i)$ is an addable node of $\lambda$
and a removable node of $\mu$. Note that $|\mu|=|\lambda|+1$.

The Young diagram $Y(\lambda)$ for a partition $\lambda=(\lambda_1,
\lambda_2, \cdots)$ is a collection of boxes arranged in
left-justified rows with $\lambda_i$ boxes in the $i$-th row of
$Y(\lambda)$. A $\lambda$-tableau $\ss$ is obtained by inserting $i,
1\le i\le n$ into $Y(\lambda)$ without repetition.  $\mathfrak S_n$
acts on $\ss$ by permuting its entries. Let $\ts^\lambda$ be the
$\lambda$-tableau obtained from $Y(\lambda)$ by adding $1, 2,
\cdots, n$ from left to right along the  rows of $Y(\lambda)$. If
$\ts^\lambda w=\ss$, write $w=d(\ss)$. Note that $d(\ss)$ is
uniquely determined by $\ss$.

A $\lambda$-tableau $\ss$ is standard if the entries in $\ss$
increase  both from left to right in each row and from top to bottom
in each column. Let $\Std_n(\lambda)$ be the set of all standard
$\lambda$-tableaux.

Given an  $\ss\in \Std_n(\lambda)$, let $\ss\!\!\downarrow_i$ be
obtained from $\ss$ by removing all  boxes containing  the entries $j$ in $\ss$ with
$j>i$. 
Let $\ss_i$ denote the partition of $i$ such that $\ss\!\!\downarrow_i$ is an $\ss_i$-tableau. 
Then $\ss=(\ss_0, \ss_1, \cdots,\ss_n)$ is a sequence of partitions such that 
$\ss_i\rightarrow\ss_{i+1}$. Conversely, if we insert $i$ into 
$\ss_i\setminus\ss_{i-1}$, we will obtain  an $\ss\in \Std_n(\lambda)$. Therefore,
there is a bijection between $\Std_n(\lambda)$ and the set of all
$(\ss_0,\ss_1, \cdots, \ss_n)$ with $\ss_0=\varnothing$,
$\ss_n=\lambda$ and  $\ss_i\rightarrow \ss_{i+1}$, $0\le i\le n-1$.

Assume that  $0\le f\le \lfloor\frac n2 \rfloor$. Let $\mathfrak
S_{n-2f}$ be the symmetric group on $\{2f+1, 2f+2, \cdots, n\}$.
Following \cite{Enyang}, let $\mathfrak B_f$ be the subgroup of
$\mathfrak S_n$ generated by $\tilde s_{i}, \tilde s_0$, where
$\tilde s_i=s_{2i} s_{2i-1}s_{2i+1}s_{2i} $, $1\le i\le f-1$, and
$\tilde s_0=s_1$. In \cite{Enyang}, Enyang  proved that
 $\mathcal D_{f, n}$  is  a  complete set of right coset
representatives of $\mathfrak B_f\times \mathfrak S_{n-2f}$ in
$\mathfrak S_n$, where
$$\mathcal D_{f, n}=\Set[90]w\in \Sym_n| $(2i+1)w<(2j+1)w$, $(2i+1)w<(2i+2)w$ for all $ 0\le
i<j<f$, and $ (k)w<(k+1) w$ for  $2f< k< n$|.
$$
For $\lambda\vdash n-2f$, let $\mathfrak S_{\lambda}$ be the Young
subgroup of $\mathfrak S_{n-2f}$ generated by $s_j$, $2f+1\le j\le
n-1$ and $j\neq 2f+\sum_{k=1}^i \lambda_k$ for all possible $i$. A
standard $\lambda$-tableau $\hat \ss$ is obtained by using $2f+i,
1\le i\le n-2f$ instead of $i$ in the usual standard
$\lambda$-tableau $\ss$. Define $d(\hat \ss)\in \mathfrak S_{n-2f}$
by declaring that $\hat \ss=\hat{\ts}^\lambda d(\hat \ss)$. If there
is no confusion, we will use $\ss$ instead of $\hat\ss$.

In \cite{GL}, Graham and Lehrer  proved that $\cba{n}$ is a cellular
algebra over a commutative ring $R$. In the current paper, we are
going to recall Enyang's cellular basis for $\cba{n}$.
We start by recalling the notion  of cellular algebra.

\begin{Defn}[\cite{GL}]\label{GL}
    Let $R$ be a commutative ring and $A$ an $R$--algebra.
    Fix a partially ordered set $\Lambda=(\Lambda,\gedom)$ and for each
    $\lambda\in\Lambda$ let $T(\lambda)$ be a finite set. Finally,
    fix $C^\lambda_{\s\t}\in A$ for all
    $\lambda\in\Lambda$ and $\s,\t\in T(\lambda)$.

    Then the triple $(\Lambda,T,C)$ is a {cell datum} for $A$ if:
    \begin{enumerate}
    \item $\set{C^\lambda_{\s\t}|\lambda\in\Lambda\text{ and }\s,\t\in
        T(\lambda)}$ is an $R$--basis for $A$;
    \item the $R$--linear map $*\map AA$ determined by
        $(C^\lambda_{\s\t})^*=C^\lambda_{\t\s}$, for all
        $\lambda\in\Lambda$ and all $\s,\t\in T(\lambda)$ is an
        anti--isomorphism of $A$;
    \item for all $\lambda\in\Lambda$, $\s\in T(\lambda)$ and $a\in A$
        there exist scalars $r_{\s\u}(a)\in R$, which are independent of $\t$,  such that
        $$aC^\lambda_{\s\t}
            =\sum_{\u\in T(\lambda)}r_{\s\u}(a)C^\lambda_{\u\t}
                     \pmod{A^{\gdom\lambda}},$$
            where
    $A^{\gdom\lambda}=R\text{--span}%
      \set{C^\mu_{\u\v}|\mu\gdom\lambda\text{ and }\u,\v\in T(\mu)}$.
    \end{enumerate}
    \noindent An algebra $A$ is a {cellular algebra} if it has
    a cell datum and in this case we call
    $\set{C^\lambda_{\s\t}|\s,\t\in T(\lambda), \lambda\in\Lambda}$
    a cellular basis of $A$.  By (c),
$A^{\vartriangleright\lambda}$ is a two-sided ideal of $A$.
\end{Defn}

Let  $\Lambda_n=\set{(f, \lambda)\mid \lambda\vdash n-2f, 0\le f\le
\lfloor\frac n 2\rfloor}$. Given  $(k, \lambda),  (f, \mu)\in
\Lambda_n$, define $(k, \lambda)\trianglelefteq (f, \mu)$  if either
$k<f$ or $k=f$ and $\lambda\trianglelefteq\mu$.
 Write
$(k, \lambda)\triangleleft(f, \mu)$, if  $(k,
\lambda)\trianglelefteq (f, \mu)$ and $(k, \lambda)\ne (f, \mu)$.

Let $I(f, \lambda)=\Std_n(\lambda)\times \mathcal D_{f, n}$. Write
$m_\lambda=e^f x_\lambda$, where $e^f=e_{1}e_{3}\cdots e_{2f-1}$ and
$x_\lambda=\sum_{w\in \mathfrak S_\lambda} w$. Denote by
$$
 C_{(\ss, u) (\ts, v)}^{(f, \lambda)} =u^{-1} d(\ss)^{-1} \mathfrak m_\lambda
d(\ts) v, \text{ for all $(\ss, u), (\ts, v)\in I(f, \lambda)$.}
$$

\begin{Theorem}\cite{Enyang1} \label{cell} Let $\cba{n}$ be a Brauer
algebra over a commutative ring $R$. Let $\sigma :
\cba{n}\rightarrow \cba{n}$ be the $R$-linear involution which
fixes $s_i, e_j$ for all $1\le i, j\le n-1$. Then
\begin{enumerate}
\item $\mathscr C_n =\{ C_{(\ss, u) (\ts, v)}^{(f, \lambda)} \mid (\ss, u), (\ts, v)\in
I(f, \lambda), \lambda\vdash n-2f,  0\le f\le \lfloor\frac
n2\rfloor\}$ is a free $R$--basis of $\cba{n}$.
\item $\sigma (C_{(\ss, u) (\ts, v)}^{(f, \lambda)})=C_{(\ts, v) (\ss,
u)}^{(f, \lambda)}$ for all $(\ss, u), (\ts, v)\in I(f, \lambda)$
and all $(f, \lambda)\in \Lambda_n$.
\item For all  $h\in \cba{n}$, and all $(\ss, u), (\ts, v)\in I(f,
\lambda)$ with $(f, \lambda)\in \Lambda_n$,
$$C_{(\ss, u) (\ts, v)}^{(f, \lambda)} h \equiv \sum_{(\us, w)\in I(f, \lambda)} a_{\us,
w} C_{(\ss, u) (\us, w)}^{(f, \lambda)} \pmod{
\cba{n}^{\vartriangleright (f, \lambda)}},$$
 where
$\cba{n}^{\vartriangleright (f, \lambda)}$ is the free $R$-submodule
generated by $ C_{(\tilde \s, \tilde u) (\tilde \ts, \tilde v)}^{(k,
\mu)}$ with $(k, \mu) \vartriangleright (f, \lambda)$ and $(\tilde
\ss, \tilde u), (\tilde \ts, \tilde v)\in I(k, \mu)$.
 Moreover, each coefficient $a_{\us, w}$ is independent of $(\ss,
u)$.
\end{enumerate}
\end{Theorem}

\textsf{From here on, all modules considered in this paper  are
right modules}.

For each $\lambda\vdash n-2f$, Enyang considered the right module
$$S^\lambda=R\text{--span}\set{\m_\lambda d(\ts) v \pmod
{\cba{n}^{\vartriangleright (f, \lambda)}}\mid (\ts, v)\in I(f,
\lambda)}$$
By the definition of a cell module in \cite{GL}, $S^{\lambda}$ is
isomorphic to the cell module $\Delta(f, \lambda)$   with respect to
the cellular basis $\mathscr C_n$, provided by Theorem~\ref{cell}.
We will identify $\Delta(f, \lambda)$ with $S^{\lambda}$. We remark
that the action of $\cba{n}$ on $\Delta(f, \lambda)$ is given by
Theorem~\ref{cell}(c).

Suppose that  $(f, \lambda)\in \Lambda_n$. An {$n$--updown
$\lambda$--tableau}, or more simply an updown $\lambda$--tableau, is
a sequence $\t=(\t_0, \t_1,\t_2,\dots,\t_n)$ of partitions, where
$\t_n=\lambda$, $\t_0=\varnothing$,  and the partition $\t_i$ is
obtained from $\t_{i-1}$ by either \textit{adding} or
\textit{removing} a box, for $i=1,\dots,n$. Let $\UPD_n(\lambda)$ be
the set of all $n$-updown $\lambda$--tableaux. If $\lambda\vdash n$,
an $n$-updown tableau becomes an $n$-up (i.e. there are no downs)
tableau. In this case, there is a bijection between
$\UPD_n(\lambda)$ and $\Std_n(\lambda)$.

For any $1\le i, j\le n$, let
$$s_{i,
j}=\begin{cases} s_{i}s_{i+1}\cdots s_{j-1}, & \text {if $j>i$,}\\
s_{i-1}s_{i-2}\cdots s_{j}, & \text{if $j<i$,}\\
1, & \text{if $i=j$.}\\
\end{cases}
$$

\begin{Defn}\label{def-m} Given  $\t\in \UPD_n(\lambda)$ with
$(f, \lambda) \in \Lambda_n $, define $f_j\in \mathbb N$ by
declaring that $\t_j\vdash j-2f_j$. Let $\mu^{(j)}=\t_j$.
Following \cite{Enyang}, define $\m_\t=\m_{\t_n}$ inductively by
declaring that
\begin{itemize}\item[(1)]
$\m_{\t_1}=1$, \item[(2)]  $ \m_{\t_i}=\sum_{j=a_{k-1}+1}^{a_k}
s_{j, i} \m_{\t_{i-1}}$ if $\t_i=\t_{i-1}\cup p$ with $p=(k,
\mu^{(i)}_k)$, and $a_l=2f_i+\sum_{j=1}^l \mu^{(i)}_j$ \item[(3)]
$\m_{\t_i}= e_{2f_i-1} s_{2f_i, i} s_{2f_i-1,  b_{k}}
\m_{\t_{i-1}} $ if $\t_{i-1}=\t_{i}\cup p$ with $p=(k,
\mu_k^{(i-1)})$,  and $b_k=2(f_i-1)+\sum_{j=1}^k
\mu^{(i-1)}_j$.\end{itemize}
\end{Defn}
 In \cite{Enyang}, Enyang  showed
  that $\m_{\t}=\m_\lambda b_{\t}$ for some $b_\t\in R \mathfrak S_n$.
The following  recursive formulae describe explicitly $b_{\t}$. Note
that  $b_\t=b_{\t_n}$ and $\t_{n-1}=\mu$.
\begin{equation}\label{des-b}
b_{\t_{n}}=\begin{cases} s_{a_k, n}b_{\t_{n-1}}, & \text{ if
$\t_{n}= \t_{n-1}\cup\{ (k, \lambda_k)\}$}\\
s_{2f, n}\sum_{j=b_{k-1}+1}^{b_k} s_{2f-1, j} b_{\t_{n-1}}, &
\text{ if $\t_{n-1}= \t_{n}\cup\{ (k, \mu_k)\}$}.\\
\end{cases}
\end{equation}

For any $(f, \lambda)\in \Lambda_n$, define
\begin{itemize}
\item $\t_{2i-1}=(1)$ and $\t_{2i}=\varnothing$ for $1\le i\le f$,
\item   $\t_{i}$ is obtained from $\hat
\ts^\lambda$ by removing the entries $j$ with $j>i$ under the
assumption  $2f+1\le i\le n$.
\end{itemize}
Then $\t=(\t_0, \t_1, \cdots, \t_n)\in \UPD_n(\lambda)$.  In this
case, we denote $\t$ by $\t^\lambda$. By definition,
$\m_{\t}=\m_{\lambda}=e^f x_\lambda$. We remark that we use  $\t$
(resp. $\ts$) to denote an updown (resp.  a standard) tableau.

Following \cite{Enyang},  let $(k_i, \mu^{(i)})\in \Lambda_{n-1},
1\le i\le m$ be such that
\begin{enumerate} \item $\mu^{(i)}\rightarrow \lambda$ if $k_i=f$,
\item $\lambda \rightarrow \mu^{(i)}$ if $k_i=f-1$,
\item $(k_1,\mu^{(1)})\vartriangleright (k_2, \mu^{(2)})\vartriangleright   \cdots
\vartriangleright (k_m, \mu^{(m)})$.
\end{enumerate}

Given $\s\in \UPD_n(\lambda)$, we
identify $\s_i$ with $(f_i, \mu^{(i)})$ if $\s_i=\mu^{(i)}\vdash
i-2f_i$. Define the partial order $\trianglelefteq $ on
$\UPD_n(\lambda)$ by declaring that $\s\trianglelefteq\t$ if
$\s_i\trianglelefteq\t_i$ for all $1\le i\le n$. Write
$\s\triangleleft\t$ if $\s\trianglelefteq\t$ and $\s\neq \t$.

Define
$$
\begin{aligned} N^{\trianglerighteq \mu^{(i)}}& =R\text{--span} \{\m_\t\pmod
{\cba{n}^{\vartriangleright (f, \lambda)}}\mid \t\in
\UPD_n(\lambda), \t_{n-1}\trianglerighteq \mu^{(i)}\},\\
N^{\triangleright \mu^{(i)}}& =R\text{--span} \{\m_\t\pmod
{\cba{n}^{\vartriangleright (f, \lambda)}}\mid \t\in
\UPD_n(\lambda), \t_{n-1}\triangleright \mu^{(i)}\}.\\
\end{aligned}
$$

Under the previous identification,  that $\t_{n-1}\trianglerighteq \mu^{(i)}$ is 
equivalent to 
either $|\t_{n-1}|<|\mu^{(i)}|$ or   $|\t_{n-1}|=|\mu^{(i)}|$  and $\t_{n-1}\trianglerighteq \mu^{(i)}$ 
under the usual dominance order.

It is well known that  $\cba{n}$ can be defined via Brauer diagrams.
In \cite{MW}, Morton and  Wassermann  proved that a Brauer algebra
defined by Brauer diagrams is isomorphic to $\cba{n}$ defined in
section~2. Therefore, the subalgebra of $\cba{n}$ generated by $s_i,
e_i, 1\le i\le n-2$, is isomorphic to $\cba{n-1}$.

\begin{Theorem}\cite{Enyang}\label{ud} Let $\cba{n}$ be a Brauer algebra
over a commutative ring $R$. Assume that $(f, \lambda)\in
\Lambda_n$.
 \begin{enumerate}
 \item
$\{\m_\t\pmod {\cba{n}^{\vartriangleright (f, \lambda)}}\mid \t\in
\UPD_n(\lambda)\}$ is an $R$-basis of $\Delta(f, \lambda)$.
\item Both $N^{\trianglerighteq \mu^{(i)}}$ and  $N^{\triangleright
\mu^{(i)}}$ are $\cba{n-1}$-submodules of $\Delta(f, \lambda)$.
\item The $R$-linear map $\phi: N^{\trianglerighteq \mu^{(i)}}/N^{\triangleright
\mu^{(i)}} \rightarrow  \Delta(k_i, \mu^{(i)})$ sending
$\m_{\t}\pmod { N^{\triangleright \mu^{(i)}}} $ to $
\m_{\t_{n-1}}\pmod{ \cba{n-1}^{\vartriangleright (k_i,
\mu^{(i)})}}$ is an isomorphism of $\cba{n-1}$-modules.
\end{enumerate}\end{Theorem}

\begin{Defn} \label{udbasis} Given  $\s, \t\in \UPD_n(\lambda)$, define $\m_{\s,
\t}=\sigma (b_{\s}) \m_{\lambda}  b_{\t}$,  where $\sigma:
\cba{n}\rightarrow \cba{n} $ is the $R$-linear anti-involution on
$\cba{n}$ given in Theorem~\ref{cell}.
\end{Defn}

\begin{Theorem}\label{Mur}  Let $\cba{n}$ be  a Brauer algebra
over a commutative ring $R$. Then
\begin{enumerate}
\item $\mathscr M_n =\{ \m_{\s, \t}\mid \s,  \t\in \UPD_n(\lambda), \lambda
\vdash n-2f, 0\le f\le \floor{ n2}\} $ is a free $R$--basis of
$\cba{n}$.
\item $\sigma (\m_{\s, \t})=\m_{\t, \s}$ for all $\s, \t\in \UPD_n(\lambda)$
and all $(f, \lambda)\in \Lambda_n$.
\item  Let $\widetilde {\cba{n}}^{\vartriangleright (f, \lambda)}$ be the free
$R$-submodule of $\cba{n}$ generated by $ \m_{\tilde \s, \tilde
\t}$ with $\tilde \s, \tilde \t\in \UPD_n(\mu)$ and $(\frac
{n-|\mu|} 2, \mu) \vartriangleright (f, \lambda)$. Then
$\widetilde {\cba{n}}^{\vartriangleright (f,
\lambda)}=\cba{n}^{\vartriangleright (f, \lambda)}$.
\item For all $\s, \t\in \UPD_n(\lambda)$ with $(f, \lambda)\in \Lambda_n$, and all $h\in \cba{n}$,
there exist scalars $a_\u\in R$ which are independent of $\s$,
such that
$$\m_{\s, \t} h \equiv \sum_{\u} a_{\u} \m_{\s,\u} \pmod{
 {\cba{n}}^{\vartriangleright (f, \lambda)}}.$$
\end{enumerate}
\end{Theorem}

\begin{proof} If $(f, \lambda)\in \Lambda_n$ is maximal, then
(c) holds  since both ${\cba{n}}^{\rhd (f, \lambda)}$ and
$\widetilde {\cba{n}}^{\rhd (f, \lambda)}$ are equal to  zero. In
general,  take a minimal element  $(k, \mu)\in \Lambda_n$ such that
$(k, \mu)\rhd (f, \lambda)$. By the induction assumption,
$\widetilde {\cba{n}}^{\rhd (k, \mu)}=\cba{n}^{\rhd (k, \mu)}$. By
Theorem~\ref{ud}(a) and  $h\m_\mu=\sigma (\m_\mu \sigma(h))$,
 $ C_{(\ss,  u) (\ts,v)}^{(k, \mu)}$ can be
expressed as a linear combination of $\m_{\s\t}$ with $\s, \t\in
\UPD_n(\mu)$  module the two-sided ideal $\cba{n}^{\rhd (k, \mu)}$.
Therefore,  $\widetilde {\cba{n}}^{\unrhd (k,
\mu)}\supseteq\cba{n}^{\unrhd (k, \mu)}$, and $\widetilde
{\cba{n}}^{\rhd (f, \lambda)}\supseteq \cba{n}^{\rhd (f, \lambda)}$.
By the similar arguments, we can verify the inverse inclusion. This
proves (c). (d) follows from Theorem~\ref{ud}(a).

Let $N$ be the $R$-module  generated by $\mathscr M_n$. By (c) and
Theorem~\ref{ud}(a),  $\m_{\s\t} h=\sigma(b_\s)( \m_\lambda b_\t
h)\in N$ for any $h\in \cba{n}$. Therefore,    $N$ is a right
$\cba{n}$-module. Since $1=x_\lambda\in N$ for $\lambda=(1,
1,\cdots, 1)\vdash n$,    $N=\cba{n}$.  Notice that $\# \mathscr M_n
=\sum_{(f, \lambda)\in \Lambda_n} \# \UPD_n(\lambda)^2=\text{rank}
\cba{n}$, $\mathscr M_n$ has to be an $R$-basis of $\cba{n}$.
Finally, (b) follows from the equality
$\sigma(\m_\lambda)=\m_\lambda$.
\end{proof}

We call $\mathscr M_n$ the {Jucys-Murphy} basis of $\cba{n}$. It is
a cellular basis of $\cba{n}$ over $R$. In order to simplify the
notation, we use $\m_\t$ instead of $\m_\t\pmod
{\cba{n}^{\vartriangleright (f, \lambda)}}$ if there is no
confusion.

In \cite{GL}, Graham and Lehrer proved that there is a symmetric
invariant bilinear form $\langle\quad, \quad \rangle: \Delta(f,
\lambda)\times  \Delta(f, \lambda)\rightarrow R$. In our case, we
use $\mathscr M_n$ to define such a bilinear form on $\Delta(f,
\lambda)$. More explicitly, $\langle \m_\s, \m_\t \rangle \in R $ is
determined by
$$\m_{\tilde\s\s}\m_{\t\tilde\t}\equiv\langle \m_\s,  \m_\t   \rangle
\m_{\tilde\s\tilde\t}\pmod { \cba{n}^{\vartriangleright (f,
\lambda)}}, \quad \text{for some $\tilde \s, \tilde\t \in
\UPD_n(\lambda)$}.$$ By  Theorem~\ref{Mur}(d), the above symmetric
invariant bilinear form is independent of $\tilde\s, \tilde\t\in
\UPD_n(\lambda)$. The Gram matrix $G_{f, \lambda}$ associated to
$\Delta(f, \lambda)$ is the $k\times k$ matrix with $k=\rank
\Delta(f, \lambda)$ such that the $(\t, \s)$-entry is $\langle
\m_\t, \m_\s\rangle$.

For any partition $\lambda=(\lambda_1, \lambda_2, \cdots, )$, let
$\lambda'=(\lambda_1', \lambda_2', \cdots, )$ be its dual partition.
Then
$$
\rank \Delta(f, \lambda)=\frac{n!(2f-1)!!}{(2f)! \prod_{(i, j)\in
\lambda} h_{i, j}^\lambda},$$  where  $h_{i,
j}^\lambda=\lambda_i+\lambda_j'-i-j+1$. In general, $ \rank
\Delta(f, \lambda)$ is a very large integer. Therefore, it is very
difficult to compute the Gram matrix $G_{f, \lambda}$ directly.


 In \cite{GL}, Graham
and Lehrer proved that a cellular algebra is (split) semisimple if
and only if the  Gram determinant associated to each cell module is
not equal to zero. Via it, we proved the following theorem.

\begin{Theorem} \cite{Rui:ssbrauer, RS:ssbrauer}\label{semisimple}
Let $\cba{n}$ be a Brauer algebra over an arbitrary field $F$.
Define $e=+\infty$ (resp. $p$) if the characteristic of $F$ is zero
(resp. $p>0$).
\begin{enumerate}
\item  Suppose $\delta\neq 0$. Then $\cba{n}$
is semisimple if and only if $$\delta\not\in \{i, -2i\mid 1\le
i\le n-2, i\in \mathbb Z\}\cup \{j\in \mathbb Z\mid 4-n\le j\le
-1\}\text{ and } e\nmid n!.$$
\item $\ccba{n}$ is semisimple if and only if $n\in \{1, 3, 5\}$
and $e\nmid n!$.
\end{enumerate}
\end{Theorem}

We remark that Theorem~\ref{semisimple} has been generalised to
cyclotomic Brauer algebras in \cite{RY, RX}.

For any  $\t\in \UPD_n(\lambda)$ with $(f, \lambda)\in \Lambda_n$,
define the residue of $k$ in $\t$ to be the scalar $c_\t(k)\in R$
such that
$$c_\t(k)= \begin{cases} \label{content}
   \frac {\delta-1} 2 + j-i,  &\text{if }\t_k=\t_{k-1}\cup \set{(i,j)},\\
    -\frac{\delta-1} 2+ i-j , &\text{if }\t_{k-1}=\t_{k}\cup\set{(i,j)}.
\end{cases}$$
We also define $$ c_\lambda(p)=\begin{cases}  \frac {\delta-1} 2 +
j-i, & \text{if $p=(i, j)$ is an addable node of $\lambda$,}\\
-\frac {\delta-1} 2 - j+i, & \text{if $p=(i, j)$ is a removable node
of $\lambda$.}\\
\end{cases}
$$
In \cite{Nazarov:brauer}, Nazarov defined the Jucys-Murphy elements
$x_i$ of $\cba{n}$ by declaring that
$$
 x_i=\begin{cases} \frac{\delta-1} 2, & \text{if $i=1,$}\\
\frac{\delta-1} 2 +\sum_{k=1}^{i-1} (s_{k,i}s_{i-1, k}- s_{k, i-2}
s_{i-1, k} e_{i-1} s_{k,i-1}s_{i-2, k}), &\text{if $2\le i\le n.$} \\
\end{cases}
$$
\begin{Remark} In the definitions of $c_\t(k)$ and $x_k$,  we assume $\frac 12\in R$.
In fact, we do not need this assumption since  $x_k$ can be defined
by using $0$ instead of $\frac{\delta-1}2$. In this situation, we
have to define $c_\t(k)=j-i$ if $\t_k= \t_{k-1}\cup \{(i, j)\}$ and
$c_\t(k)=1-\delta-j+i$ if $\t_{k-1}= \t_{k}\cup \{(i, j)\}$.
\end{Remark}

Nazarov  proved the following Lemma for $\cba{n}$ over $\mathbb C$. It holds  in general \cite{NW}.
\begin{Lemma}\cite[\S4]{Nazarov:brauer}\label{degen}
 {\small \begin{multicols}{2}
\begin{enumerate}
\item $e_ix_j=x_je_i$,  $j\ne i,i+1$,
\item $s_ix_j=x_js_i$,  $j\ne i,i+1$,
\item $s_ix_i-x_{i+1}s_i=e_i-1$,  $1\le i<n$,
\item $x_is_i-s_ix_{i+1}=e_i-1$,  $1\le i<n$,
\item $e_i(x_i+x_{i+1})=0$, $1\le i<n$,
\item $(x_i+x_{i+1})e_i=0$, $1\le i<n$,
\item $x_ix_j=x_jx_i, 1\le i, j\le n$.
\end{enumerate}
\end{multicols}}
\end{Lemma}

\begin{Theorem}\cite[10.7]{Enyang} \label{xproduct} Given $\s, \t\in \UPD_n(\lambda)$,
 with $(f, \lambda)\in \Lambda_n$, let $\m_{\s\t}\in \cba{n} $ be defined in Definition~\ref{udbasis}.
 Then  $$\m_{\s\t} x_k\equiv c_\t(k)\m_{\s\t} +
\sum_{\substack{\u\in \UPD_n(\lambda)\\ \u\triangleright\t\\}} a_\u
\m_{\s\u} \pmod {\cba{n}^{\vartriangleright (f, \lambda)}}
.$$\end{Theorem}

Theorem~\ref{xproduct} was proved in \cite[10.7]{Enyang} under the
assumption $\s=\t^\lambda$. In general, it follows from this special
case
  since multiplying an element on the left
 side of $\m_\t\pmod {\cba{n}^{\vartriangleright (f, \lambda)}}$
  gives a homomorphism of right $\cba{n}$-module.

\section{Orthogonal representations for $\cba{n}$}
We always assume that $F$ is a field which satisfies the
assumptions~\ref{assumption} in this section. The main purpose of
this section is to construct an orthogonal basis. We remark that
many results in this section are  motivated by \cite{M:gendeg}.

\begin{Assumption}\label{assumption} Suppose that  $F$ is a field of
characteristic $p$ where  either $p=0$ or $p>2n$. Suppose
$\delta\in F$ such that $|c|\ge 2n+1$ whenever $\delta-c\cdot 1_F=0$
for some $c\in \mathbb Z$.\end{Assumption}

For example, $\mathbb C(\delta)$  satisfies
assumption~\ref{assumption}, where $\delta$ is an indeterminate.

\begin{Defn}\label{tilde} Suppose $1\le k\le n-1$ and
$(f, \lambda)\in \Lambda_n$.
Define an equivalence relation $\simk$ on $\UPD_n(\lambda)$ by
declaring that $\t\simk \s$ if $\t_j=\s_j$ whenever $1\le j\le n$
and $j\neq k$, for $\s,\t\in\UPD_n(\lambda)$.
\end{Defn}

The following result is a special case of \cite[4.2]{NW}.

\begin{Lemma} Suppose $\s\in \UPD_n(\lambda)$ with
$\s_{k-1}=\s_{k+1}$. Then there is a bijection between the set of
all addable and removable nodes of $\s_{k+1}$ and the set
$\set{\t\in \UPD_n(\lambda)\mid \t\simk \s}$.
\end{Lemma}

Suppose $\lambda$ and $\mu$ are partitions. We write
$\lambda\ominus\mu=\alpha$ if either $\lambda\supset \mu$ and
$\lambda\setminus \mu=\alpha$ or $\lambda\subset \mu$ and
$\mu\setminus \lambda=\alpha$.

\begin{Lemma} \label{cont} Suppose that $F$ is a field which satisfies   the assumption~\ref{assumption}.
Assume  that $\s, \t\in \UPD_n(\lambda)$ with $(f, \lambda)\in
\Lambda_n$.
\begin{enumerate} \item $\s=\t$ if and only if $c_\s(k)=c_\t(k)$
for $1\le k\le n$.
\item If  $\t_{k-1}\neq \t_{k+1}$, then   $c_\t(k) \pm c_\t(k+1)\neq 0$.
\item If $\t_{k-1}=\t_{k+1}$, then $c_\t(k)\pm c_\s(k)\neq 0$
whenever $\s\simk\t$ and $\s\neq \t$.
\item  If $\t_{k-1}=\t_{k+1}$, then $2c_\t(k)+1\neq 0$.
 \end{enumerate}
\end{Lemma}

\begin{proof} The ``only if " part of (a) is trivial. We  prove the ``if " part of (a) as follows.
(b)-(d) can be proved similarly. We  leave the details to the
reader.

By the definition of the $n$-updown tableau, $\t_1=(1)$ for any
$\t\in \UPD_n(\lambda)$. Assume $\t_{k-1}=\s_{k-1}$ and
$\t_k{\ominus}\t_{k-1}=(i,j)$ and $\s_k{\ominus}\s_{k-1}=(i',j')$.

If the sign of $\delta$ in $c_\s(k)$ and $c_\t(k)$ are different,
then, $n\ge 3$ and $(i, j)\neq (i', j')$ since one node of a Young
diagram can not be an addable and a removable node simultaneously.
We  write $ c_\t(k)= \frac {\delta-1} 2 + j-i $ and $ c_\s(k)=
-\frac{\delta-1} 2+ i'-j'$ without loss of generality. In this
situation, $\lambda\vdash n-2f$ and  $f>0$. The maximal (resp.
minimal) value of $(i-j)+(i'-j')$ is $2n-5$ (resp. $5-2n$). In the
first case, $\lambda=(1,1, \cdots, 1)\vdash n-2$ and $\{(i, j), (i',
j')\}= \{(n-1, 1), (n-2, 1)\}$. In the second case,
$\lambda=(n-2)\vdash n-2$ and  $\{(i, j), (i', j')\}=\{(1, n-1), (1,
n-2)\}$. In any case, $\delta=1+i'-j'+i-j$ if $c_\s(k)=c_\t(k)$. We
have $|1+(i'-j'+i-j)|\le 2n-4$, which contradicts  the
assumption~\ref{assumption}.

Suppose that the sign of $\delta$ in $c_\s(k)$ and $c_\t(k)$ are the
same. Then both $(i, j)$ and $(i', j')$ are either  addable nodes or
 removable nodes of $\t_{k-1}=\s_{k-1}$. If $\s_k\neq \t_k$, then
$(i, j)\neq (i', j')$. In this situation, such two nodes can not be
in the same diagonal of the partition  which is obtained from
$\t_{k-1}$ by adding (if $\t_k\supset \t_{k-1}$) or removing (if
$\t_k\subset \t_{k-1}$) $(i, j)$ and $(i', j')$. Therefore, $i-j\neq
i'-j'$ and $c_\s(k)\neq c_\t(k)$, a contradiction.

\end{proof}

Following \cite{Mur}, we make the following definition.

\begin{Defn} \label{Felement} Let  $\cba{n}$ be  a Brauer
algebra over a field $F$. Define
\begin{enumerate}
\item  $\mathscr R(k)=\set{ c_\t(k)\mid \t\in \UPD_n(\lambda), (f,
\lambda)\in \Lambda_n}$ for any $k\in \mathbb Z$, $1\le k\le n$.
\item
$F_\t =\prod_{k=1}^n \prod_{\substack{r\in \mathscr R(k)\\
c_\t(k)\neq r}} \frac {x_k-r} {c_\t(k)-r} $
 \item
 $
f_{\s\t}=F_\s \m_{\s\t} F_\t$,
 \item
$f_\s= \m_\s  F_\s$, \end{enumerate} where  $\s, \t\in
\UPD_n(\lambda)$ and $(f, \lambda)\in \Lambda_n$.\end{Defn}

In (d), we use $\m_\s$ instead of $\m_\s\pmod
{\cba{n}^{\vartriangleright (f, \lambda)}}$. Therefore, (d) should
be read as $f_\s\equiv  \m_\s  F_\s  \pmod
{\cba{n}^{\vartriangleright (f, \lambda)}}$.  In what follows, we
will omit $\pmod {\cba{n}^{\vartriangleright (f, \lambda)}}$ for the
simplification in exposition and notation.

\begin{Lemma}\label{fxproduct} Let $\cba{n}$ be a Brauer
algebra over a field $F$.
 Suppose that  $\t\in \UPD_n(\lambda)$ with  $(f, \lambda)\in \Lambda_n$.
\begin{enumerate}\item  $f_\t=\m_\t + \sum_{\s\in \UPD_n(\lambda)} a_\s
\m_\s$, and $\s\vartriangleright \t$ if $a_\s\neq 0$.
 \item $\m_\t=f_\t + \sum_{\s\in \UPD_n(\lambda)} b_\s
f_\s$, and $\s\vartriangleright\t$ if $b_\s\neq 0$.
\item $f_\t x_k=c_\t(k) f_\t$, for any integer $k$, $1\le k\le n$.
\item $f_\t  F_\s =\delta_{\s\t} f_\t$ for all $\s\in
\UPD_n(\mu)$ with $(\frac {n-|\mu|}{2}, \mu)\in \Lambda_n$.
\end{enumerate}
\end{Lemma}

\begin{proof} (a)  follows from  Theorem~\ref{xproduct}. Since the
transition matrix from $\{f_\s\mid  \s\in \UPD_n(\lambda)\}$ to
$\{\m_\s\mid  \s\in \UPD_n(\lambda)\}$ is upper uni-triangular, so
is its inverse. This prove (b). (c)  can be proved by the arguments
in \cite[3.35]{M:ULect}. The only difference is that we need use $\#
\UPD_n(\lambda)$ instead of $\#\Std_n(\lambda)$ in the proof of
\cite[3.35]{M:ULect}.  (d) follows from (c).
\end{proof}

By Lemma~\ref{fxproduct}(a), $\set{f_\t\mid \t\in \UPD_n(\lambda)}$
is an $F$-basis for $\Delta(f, \lambda)$. Such a basis will be
called an orthogonal basis of $\Delta(f, \lambda)$.

The following two results follow from Lemma~\ref{fxproduct}.

\begin{Cor}\label{equal} Let $\cba{n}$ be a Brauer
algebra over a field $F$. Let $G_{f, \lambda}$ (resp. $\tilde G_{f,
\lambda}$) be the Gram matrix associated to the cell module
$\Delta(f, \lambda)$ with $(f, \lambda)\in \Lambda_n$, which is
defined via its Jucys-Murphy (resp. orthogonal) basis. Then $\det
G_{f, \lambda}=\det \tilde G_{f, \lambda}$.
\end{Cor}

\begin{Cor}\label{f-filtration} Let  $\cba{n}$ be a Brauer
algebra over a field $F$. Keep the setup in Theorem~\ref{ud}.
\begin{enumerate}
\item $\{f_\t\mid \t\in \UPD_n(\lambda), \t_{n-1}\trianglerighteq
\mu^{(i)}\}$ is an $F$-basis of $N^{\trianglerighteq \mu^{(i)}}$.
\item $\{f_\t\mid \t\in
\UPD_n(\lambda), \t_{n-1}\triangleright \mu^{(i)}\}$ is an
$F$-basis of $N^{\triangleright \mu^{(i)}}$.
\item The $F$-linear
map $\phi: N^{\trianglerighteq \mu^{(i)}}/N^{\triangleright
\mu^{(i)}} \rightarrow \Delta(k_i, \mu^{(i)})$ sending $f_{\t}+
N^{\triangleright \mu^{(i)}} $ to $ f_{\t_{n-1}}$ is an
isomorphism of  $\cba{n-1}$-modules.
\end{enumerate} \end{Cor}
The following result can be proved by using the arguments in
\cite{M:gendeg}.

\begin{Lemma}\label{fst} Let $\cba{n}$ be a Brauer
algebra over a field $F$.
\begin{enumerate}
\item $\set{f_{\s\t} \mid \s, \t\in \UPD_n(\lambda), (f, \lambda)\in \Lambda_n}$
 is a  basis of $\cba{n}$ over $F$.
\item $f_{\s\t} x_k=c_\t(k) f_{\s\t}$ for any $\s, \t\in
\UPD_n(\lambda)$ with $(\frac {n-|\lambda|}{2}, \lambda)\in
\Lambda_n$.
\item $f_{\s\t} F_\u=\delta_{\t\u} f_{\s\t}$  for any $\s, \t\in
\UPD_n(\lambda)$ with $(\frac {n-|\lambda|}{2}, \lambda)\in
\Lambda_n$ and any $\u\in \UPD_n(\mu)$ with $(\frac {n-|\mu|}{2},
\mu)\in \Lambda_n$.
\item Suppose  $\s, \t\in
\UPD_n(\lambda)$ with $(\frac {n-|\lambda|}{2}, \lambda)\in
\Lambda_n$ and  $\u, \v\in \UPD_n(\mu)$ with $(\frac {n-|\mu|}{2},
\mu)\in \Lambda_n$. Then  $f_{\s\t} f_{\u\v}=\delta_{\t \u}
\langle f_\t, f_\t  \rangle f_{\s\v}$ for some scalar  $\langle
f_\t, f_\t \rangle\in F$.
\end{enumerate}
\end{Lemma}

Under the assumption~\ref{assumption}, $\cba{n}$ is semisimple over
$F$ (see Theorem~\ref{semisimple}).  Therefore, $\det G_{f,
\lambda}\neq 0$ for all $(f, \lambda)\in \Lambda_n$. By
Corollary~3.7, \begin{equation}\label{fneq0}\langle f_\t,
f_\t\rangle\neq 0, \text{ for all $\t\in
\UPD_n(\lambda)$.}\end{equation}

\begin{Lemma}\label{epro} Let  $\cba{n}$ be a Brauer
algebra over a field $F$. Suppose $\t\in \UPD_n(\lambda)$ and $1\le
k\le n-1$.
\begin{enumerate}
\item Let  $f_\t s_k=\sum_{\s\in \UPD_n(\lambda)} s_{\t\s}(k)
f_\s$. If $s_{\t\s}(k)\neq 0$, then  $\s\simk \t$.
\item Let $f_\t e_k=\sum_{\s\in \UPD_n(\lambda)} e_{\t\s}(k) f_\s$.
If $e_{\t\s}(k)\neq 0$,  then $\s\simk \t$.
\end{enumerate}
\end{Lemma}
\begin{proof} By Lemma~2.10(a)-(b), $c_\s(j)=c_\t(j)$ for $j\neq k, k+1$ 
if either
$s_{\t\s}(k)\neq 0$ or $e_{\t\s}(k)\neq 0$.  Notice that
$\s_0=\t_0=\varnothing$ and $\s_n=\t_n=\lambda$. Applying
Lemma~\ref{cont}(a) to the  sequences $(\s_0, \s_1, \cdots,
\s_{k-1})$ and $(\s_{k+1}, \cdots, \s_{n})$,  we have $\s\simk \t$,
as required.
\end{proof}

\begin{Lemma} \label{ffe} Let $\cba{n}$ be  a Brauer
algebra over a field $F$. Suppose $\t\in \UPD_n(\lambda)$. If
$\t_{k-1} \neq \t_{k+1}$, then $f_\t e_k=0$.\end{Lemma}
\begin{proof}

By Lemma~\ref{degen}(f),  $f_\t (x_k+x_{k+1})
e_k=(c_\t(k)+c_\t(k+1))   f_\t e_{k} =0$. If $f_\t e_k\neq 0$, then
$c_\t(k)+ c_\t (k+1)=0$, which contradicts
Lemma~\ref{cont}(b).\end{proof}

The following result can be proved easily.

\begin{Lemma} \label{udexist}  Suppose $\t \in \UPD_n(\lambda)$
with $\t_{k-1}\neq \t_{k+1}$.
\begin{enumerate}\item  If $\t_k\ominus \t_{k-1}$ and
$\t_{k+1} \ominus \t_k$ are neither in the same row nor in the
same column, then there is a unique up-down tableau in
$\UPD_n(\lambda)$, denoted by $\t s_k $, such that  $\t s_k\simk
\t$ and $c_{\t}(k)=c_{\t s_k}(k+1)$ and $c_{\t}(k+1)=c_{\t
 s_k}(k)$.\item  If $\t_k\ominus \t_{k-1}$ and $\t_{k+1} \ominus
\t_k$ are either in the same row or in the same column, then
$\s\simk \t$ if and only if $\s=\t$.
\end{enumerate}
\end{Lemma}

\begin{Lemma}\label{fss} Let $\cba{n}$ be a Brauer
algebra over a field $F$.
 Suppose $\t\in \UPD_n(\lambda)$ with $\t_{k-1} \neq \t_{k+1}$.
If $\t s_k$ exists, then   $$f_\t s_k=\frac {1}{c_\t (k+1)-c_\t(k)}
f_\t + s_{\t,\t s_k} (k) f_{\t s_k}, $$ where
$$s_{\t,\t s_k}(k)=\begin{cases} 1, & \text{if $\t
s_k\vartriangleleft \t$},\\
 1-\frac {1}{(c_\t (k+1)-c_\t(k))^2},
&\text{if $\t s_k \vartriangleright \t$.}\\
\end{cases}
$$
\end{Lemma}

\begin{proof} Write $f_\t s_k=\sum_{\s\in \UPD_n(\lambda)}
s_{\t\s}(k) f_\s$. Since $s_k x_k=x_{k+1}s_k +e_k-1$ and $f_\t
e_k=0$ (see Lemma~\ref{ffe}), $f_\t s_k x_k=c_\t(k+1) f_\t
s_k-f_\t$.  By Lemma~\ref{epro}, $\s\simk\t$ if $s_{\t\s}(k)\neq 0$.
 Comparing the coefficient of $f_\t$ in $f_\t s_k x_k$,
we obtain $$s_{\t\t}(k)=\frac {1}{c_\t (k+1)-c_\t(k)}.$$ 
The previous formula makes sense by Lemma~3.4b. 
Suppose
$\s\neq \t$. Comparing the coefficient of $f_\s$ in $f_\t s_k x_k$,
we have
 $s_{\t\s}(k) c_{\s}(k)=s_{\t\s}(k) c_{\t}(k+1)$ forcing
 $ c_{\s}(k)= c_{\t}(k+1)$. If we compare the coefficient of
$f_\s$ in $f_\t x_ks_k$, we obtain $ c_{\s}(k+1)= c_{\t}(k)$. In
other words, $\s=\t s_k$ if $s_{\t\s}(k)\neq 0$ with $\s\neq \t$. We
claim that
\begin{enumerate} \item $\m_\t s_{k}=\m_{\t s_{k}}$,
\item $f_{\t s_k}$ does not appear in the expression of $f_\u s_k$
with $\u\vartriangleright \t$.
\end{enumerate}
 By Lemma~\ref{fxproduct}(a)-(b), and the above claims, we have
  $s_{\t\s}(k)=1$ for  $\s=\t s_k\vartriangleleft \t$.

 First, we prove (b). By Lemma~\ref{epro}(a),
 $f_\u s_k=\sum_{\v\in \UPD_n(\lambda)} s_{\u\v}(k)
f_\v$, and $\v\simk \u$ if $s_{\u\v}(k)\neq 0$.
 If $\v=\t s_k$, then $\u\simk \t$,  forcing $\u_{k-1}\neq \u_{k+1}$.
Therefore, $\u\in \{\t, \t s_k\}$ which contradicts  the assumption
$\u \vartriangleright \t\vartriangleright \t s_k$.

Now we prove (a). Note that $\m_\t s_{k}=\m_{\t s_{k}}$ if and
only if $\m_\t=\m_{\t s_{k}} s_{k}$. Therefore, (a) holds 
although $\t s_k\vartriangleright \t$.

By Theorem~\ref{ud}(c), the coefficient of $\m_{\t s_k}$ in $\m_\t
s_k$ is completely determined by $\m_{\t_{k+1}} s_k$. So, we can
assume  $k=n-1$ without loss of  generality when we prove (a).

There are four cases we have to discuss. In any case, $\t_{n}\ominus
\t_{n-1}$ and $\t_{n-1}\ominus \t_{n-2}$ are neither in the same row
nor in the same column of $\t_{n-1}$. Otherwise, $\t s_{n-1}$ does
not exist. Suppose $\t_n=\lambda$ and $\t_{n-1}=\mu$.

\Case{1. $\t_n\supset \t_{n-1}\supset\t_{n-2}$} There are some $k,
l\in \mathbb N$ such that $\t_{n}\setminus \t_{n-1}=(k, \lambda_k)$
and $\t_{n-1}\setminus \t_{n-2}=(l, \mu_l)$. If we write
$a=2f+\sum_{j=1}^k \lambda_j$ and $b=2f+\sum_{j=1}^l \mu_j$, then
$\m_{\t_n}=\m_\lambda s_{a, n} s_{b, n-1} b_{\t_{n-2}}$. Note that
$b_{\t_{n-2}}\in \mathfrak S_{n-2}$ (see (\ref{des-b})),
$b_{\t_{n-2}} s_{n-1}=s_{n-1} b_{\t_{n-2}}$. Therefore,  $$ \m_\t
s_{n-1} =\m_\lambda s_{a, n} s_{b, n-1}s_{n-1} b_{\t_{n-2}}
 =\m_\lambda s_{b, n } s_{a-1, n-1} b_{\t_{n-2}}=\m_{\t s_{n-1}}$$
  if $a>b$.
The result for  $a<b$ is still true since we can  switch $\t$ to $\t
s_{n-1}$ in the above argument. We remark that $a>b$ if  $\t
s_{n-1}\vartriangleleft \t$.

\Case{2. $\t_n\subset \t_{n-1}\supset\t_{n-2}$} Then
$\t_{n-1}\setminus \t_{n}=(k, \mu_k)$ and $\t_{n-1}\setminus
\t_{n-2}=(l, \mu_l)$.  Write $a_i=2(f-1)+\sum_{j=1}^i \mu_j$.  If
$l>k$, then $a_l>a_k>a_{k-1}$. Therefore,
$$\begin{aligned} \m_{\t_n}s_{n-1}& =\m_\lambda s_{2f, n}
\sum_{j=a_{k-1}+1}^{a_k} s_{2f-1, j} s_{a_l, n} b_{\t_{n-2}}\\ &
=\m_{\lambda} s_{a_l+1, n} s_{2f, n-1}
\sum_{j=a_{k-1}+1}^{a_k} s_{2f-1, j}b_{\t_{n-2}} =\m_{\t s_{n-1}}.\\
\end{aligned}
$$
If $l<k$, then $a_l\le a_{k-1}<a_k$. By braid relation,
$$\begin{aligned} \m_{\t_n}s_{n-1}&=\m_\lambda s_{2f, n}
\sum_{j=a_{k-1}+1}^{a_k} s_{2f-1, j} s_{a_l, n-1} s_{n-1}
 b_{\t_{n-2}}\\ &=\m_{\lambda} s_{2f, n} s_{a_l+1, n}
 \sum_{j=a_{k-1}+1}^{a_k}  s_{2f-1, j-1} b_{\t_{n-2}}\\
&=\m_\lambda  s_{a_l+2, n} s_{2f, n-1} \sum_{j=a_{k-1}}^{a_k-1}
s_{2f-1, j}b_{\t_{n-2}}=\m_{\t s_{n-1}}.\\
\end{aligned}$$
  Finally, we remark that $\t s_{n-1} \vartriangleleft\t$ is not
compatible with $\t_n\subset \t_{n-1}\supset\t_{n-2}$.

\Case{3. $\t_n\supset \t_{n-1}\subset\t_{n-2}$} The result follows
from case~2 if  we switch  $\t s_{n-1}$ to $\t$ in case 2.

\Case{4. $\t_{n}\subset\t_{n-1}\subset \t_{n-2}$} Suppose
$\t_{\n-1}= \t_n\cup (l, \nu_l) $ and $\t_{n-2}= \t_{n-1}\cup (k,
\nu_k)$ and $\t_{n-2}=\nu$. If $\t s_{n-1}\vartriangleleft \t$,
then  $k>l$. Let $b=2(f-1)+\sum_{j=1}^l \nu_j$ and $a=2(f-2)
+\sum_{j=1}^k \nu_j$. Then  $ a+1\ge b$. In this situation,
$$\begin{aligned}  \m_\t
& =e_{2f-1} s_{2f, n} s_{2f-1, b} e_{2f-3} s_{2f-2, n-1} s_{2f-3,
a}
\m_{\t_{n-2}}, \\
\m_{\t s_{n-1}}&=e_{2f-1} s_{2f, n} s_{2f-1, a+1} e_{2f-3}
s_{2f-2, n-1} s_{2f-3, b-2}\m_{\t_{n-2}}.\\ \end{aligned}
$$
On the other hand, by direct computation (in fact, one can verify
the following equalities easily if one uses Brauer diagrams in, e.g.
\cite{Wenzl:ssbrauer}),
$$\begin{aligned} & e_{2f-1}e_{2f-3}  s_{2f, n} s_{2f-1, b}  s_{2f-2, n-1} s_{2f-3,
a} s_{n-1}\\  = & e_{2f-1}e_{2f-3} s_{2f-1, n}    s_{2f-2, n-1}
s_{2f-2, a} s_{2f-3, b-2}\\ = & e_{2f-1}s_{2f, n}s_{2f-1,
a+1}e_{2f-3} s_{2f-2, n-1} s_{2f-3, b-2}.\end{aligned}
 $$
Therefore, $\m_{\t} s_{n-1}=\m_{\t s_{n-1}}$. This completes the
proof of (a).

 By considering $f_\t s_k^2=f_\t$ for $\t
s_k\vartriangleleft \t$, we obtain the result for $\t
s_k\vartriangleright \t$.
\end{proof}

\begin{Lemma} Let $\cba{n}$ be a Brauer
algebra over a field $F$. Suppose $\t\in \UPD_n(\lambda)$ with
$\t_{k-1} \neq \t_{k+1}$.
\begin{itemize}\item [(a)] If $\t_k\ominus \t_{k-1}$ and $\t_k\ominus \t_{k+1}$ are in the same
row,   then $f_\t s_k=f_\t$.
\item[(b)] If $\t_k\ominus \t_{k-1}$ and $\t_k\ominus \t_{k+1}$ are in the same column,  then  $f_\t s_k =-f_\t$.
\end{itemize}
\end{Lemma}

\begin{proof}Write $f_\t s_k=\sum_{\s\in \UPD_n(\lambda)}
s_{\t\s}(k) f_\s$. By Lemma~\ref{epro}, $\s\simk \t$ if
$s_{\t\s}(k)\neq 0$. Under the assumption,  $\s\simk \t$ if and only
if $\s=\t$. Therefore,  $f_\t s_k=s_{\t\t} (k) f_\t$, and $f_\t
=s_{\t\t} (k) f_\t s_k$. We have $s_{\t\t} (k)\in \{-1, 1\}$. Since
 $f_\t x_{k+1}=f_\t s_k x_k s_k+f_\t s_k$,  $s_{\t\t} (k)=1$
(resp. $s_{\t\t} (k)=-1$) if $\t_k\ominus \t_{k-1}$ and
 $\t_k \ominus \t_{k+1}$ are in the same row (resp. column).
\end{proof}


\begin{Lemma}\label{ek}
Let  $\cba{n}$ be a Brauer algebra over $F$ such that \cite[4.12]{NW} for $\cba{n}$  holds.
 Suppose that  $\t\in \UPD_n(\lambda)$ and
$\t_{k-1}=\t_{k+1}$. Then
\begin{enumerate}
\item $
    e_{\t\t}(k)=\big(2c_\t(k)+1\big)
             \prod_{\substack{\s\simk \t\\ \s\neq \t}
             }\frac{c_\t(k)+c_\s(k)}{c_\t(k)-c_\s(k)}\neq 0$.

\item $e_{\t\s}(k)e_{\u\u}(k)=e_{\t\u}(k)e_{\u\s}(k)$ for any
$\s, \u\in \UPD_n(\lambda)$ with $\t\simk\s\simk \u$.
\end{enumerate}
\end{Lemma}

\begin{proof} Given $(f, \lambda)\in \Lambda_n$.  We claim
 $e_{\t\t}(k)\neq 0$ for any $\t\in \UPD_n(\lambda)$
 with $\t_{k-1}= \t_{k+1}$.

In \cite[\S4]{NW}, Ariki, Mathas and the first author constructed
the seminormal representations for cyclotomic Nazarov-Wenzl algebras
$\mathscr W_{r, n}$ under the assumption \cite[4.3]{NW} and \cite[4.12]{NW}, 
the latter is called the root conditions. In
particular, $\mathscr W_{1, n}$ is $\cba{n}$. In this case,   \cite[4.3]{NW} is 
our assumption~3.1. 
Therefore, the seminormal
representations $S^{(l, \mu)}$ (denoted by $\Delta(\mu)$ in \cite{NW})
 for $\cba{n}$  were constructed for
all $(l, \mu)\in \Lambda_n$. More explicitly, each $S^{(l, \mu)}$
has  a basis $\{v_\t\mid \t\in\UPD_n(\mu)\}$ such that
\begin{itemize}\item[(1)] $v_\t x_k=c_\t(k) v_\t$.
\item [(2)] $v_\t e_k=\sum_{ \s\simk
\t} \tilde e_{\t\s}(k) v_\s$ with $\tilde e_{\t\t}(k)\neq
0$.\end{itemize} Since we are assuming that  $F$ satisfies the
assumption~\ref{assumption}, $\cba{n}$ is semisimple. We have
$Fv_\t=\cap_{k=1}^n S^{(l, \mu)}_{c_\t(k)}$, where $S^{(l,
\mu)}_{c_\t(k)}$ is the eigenspace of $S^{(l, \mu)}$ with respect to
eigenvalue  $c_\t(k)$ of $x_k$. Since $\{S^{(l, \mu)}\mid (l,
\mu)\in \Lambda_n\}$ consists of all pair-wise non-isomorphic
irreducible $\cba{n}$-modules when $\cba{n}$ is semisimple,
$\Delta(f, \lambda)\cong S^{(l, \mu)}$ for some $(l, \mu)\in
\Lambda_n$. Let $\phi$ be the corresponding isomorphism. Then
$\phi(f_\t)\in \cap_{k=1}^n S^{(l, \mu)}_{c_\t(k)}$ for some $(l,
\mu)\in \Lambda_n$. Therefore,  $\phi(f_t)$ is equal to $v_\s\in
S^{(l, \mu)} $ for some $\s\in \UPD_n(\mu)$ up to a non-zero scalar.
This implies that $e_{\t\t}(k)\neq 0$.

Suppose $k\ge1$ and that $a\ge0$. Let $Z\(\cba{k-1}\)$ be the center
of $\cba{k-1}$.  It is proved in  \cite[4.15]{NW}
 that there exist elements
$\delta^{(a)}_k$ in $Z\(\cba{k-1}\)\cap F[x_1,\dots,x_{k-1}]$ such
that $e_kx_k^ae_k=\delta^{(a)}_k e_k$. Moreover, the generating
series
    $W_k(y)=\sum_{a\ge0}\delta_k^{(a)}y^{-a}$ satisfies
    $$W_{k+1}(y)+y-\half=\( W_k(y)+y-\half\)
     \frac{(y+x_k)^2-1}{(y-x_k)^2-1}\cdot \frac{(y-x_k)^2}{(y+x_k)^2}
     .$$
Comparing the coefficient of $f_\s$ on both sides of the identity
$f_\t e_k  W_{k}(y)=f_\t e_k \frac y{y-x_k} e_k$, we have
$$
{W_k(y, \s)}{y^{-1}} e_{\t\s}(k) =\sum_{\s\simk\t\simk \u}\frac{
e_{\t\u}(k) e_{\u\s} (k) } {y-c_\u(k)}.$$ Therefore,
$$ e_{\t\s}(k)\cdot
Res_{y=c_\u(k)} { W_k(y, \s)y^{-1}} =e_{\u\s}(k) e_{\t\u} (k).
$$
Since  $e_{\t\t}(k)\neq 0$,  $e_{\t\t}(k)=Res_{y=c_\t(k)}{ W_k(y,
\t)}{y^{-1}}$ by assuming that $\t=\s=\u$.

If $\s\simk \u$, then $c_\s(j)=c_\u(j)$ for $j\le k-1$. Note that
$W_k(y)\in F[x_1, x_2, \cdots, x_{k-1}]$, ${W_k(y, \s)}=W_k(y, \u)$.
Therefore,
 $e_{\t\s}(k)e_{\u\u} (k) =e_{\u\s}(k) e_{\t\u} (k)$.
 In order to prove (a), we need to  prove
 $$Res_{y=c_\t(k)} {W_k(y,
\t)}/{y}=\big(2c_\t(k)+1\big)
             \prod_{\substack{\s\simk \t\\ \s\neq \t}
             }\frac{c_\t(k)+c_\s(k)}{c_\t(k)-c_\s(k)},$$
which is \cite[4.8]{NW} for $r=1$.
\end{proof}

It is proved in  \cite[5.4]{NW} that the real field   $\mathbb R$  
satisfies  the assumption~\ref{assumption} and \cite[4.12]{NW}. 
If we do not use Lemma~\ref{ek}(a), we can remove the assumption \cite[4.12]{NW}.

\begin{Lemma}  Let $\cba{n}$ be a Brauer
algebra over a field $F$. Suppose $\t\in\UPD_n(\lambda)$ with
$\t_{k-1}=\t_{k+1}$.
 \begin{enumerate}
\item $f_\t e_k=\sum_{\s\simk\t}
e_{\t\s} (k) f_\s$.  Furthermore, $\langle f_\s,  f_\s\rangle
e_{\t\s} (k)=\langle f_\t,  f_\t\rangle e_{\s\t} (k)$.
\item
$f_\t s_k=\sum_{ \s\simk\t} s_{\t\s} (k) f_\s$. Furthermore,
$s_{\t\s}(k)=\frac{e_{\t\s}(k)-\delta_{\t\s}}{c_\t(k)+c_\s(k)}$.
\end{enumerate}
\end{Lemma}
\begin{proof} (a) follows from $\langle f_\t e_k, f_\s\rangle=\langle f_\t, f_\s
e_k\rangle$, and $\langle f_\s, f_\t\rangle=\delta_{\s\t} \langle
f_\t, f_\t\rangle$.  (b) follows from (a) and
Lemma~2.10(d).\end{proof}

We organise the above results as follows. Such a result can be
considered as a generalisation of Dipper-James theorem for Hecke
algebras  of type $A$ \cite{DJ}.

\begin{Theorem} \label{gram} Let $\cba{n}$ be a Brauer
algebra over a field $F$. Assume that $\t\in \UPD_n(\lambda)$.

\begin{enumerate} \item Suppose
$\t_{k-1}\neq \t_{k+1}$. Then {\small $$ f_\t s_k=\begin{cases}
f_\t, &\text{\small{if $\t_k\ominus \t_{k-1}$
and $\t_k\ominus \t_{k+1}$ are in the same row  of $\t_k$,}}\\
-f_t, &\text{\small {if $\t_k\ominus \t_{k-1}$
and $\t_k\ominus \t_{k+1}$ are in the same column  of $\t_k$,}}\\
 \frac {1}{c_\t
(k+1)-c_\t(k)} f_\t + f_{\t s_k}, & \text{\small if $\t s_k\in \UPD_n(\lambda)$ and $\t s_k\vartriangleleft\t, $}\\
\frac {1}{c_\t (k+1)-c_\t(k)} f_\t + c f_{\t s_k}, & \text{\small
if $\t s_k\in \UPD_n(\lambda)$ and $\t s_k\vartriangleright \t,$}\\
\end{cases}
$$}
where $c=\frac {(c_\t (k+1)-c_\t(k)+1)( c_\t
(k+1)-c_\t(k)-1)}{(c_\t (k+1)-c_\t(k))^2}$.
\item Suppose  $\t_{k-1}\neq \t_{k+1}$.
Then  $f_\t e_k=0$.
\item Suppose  $\t_{k-1}=\t_{k+1}$. Then
$f_\t e_k=\sum_{\s\simk\t }e_{\t\s} (k) f_\s$ with
  $$\langle f_\s,  f_\s\rangle e_{\t\s}
(k)=\langle f_\t , f_\t\rangle e_{\s\t} (k).$$
\item Suppose  $\t_{k-1}=\t_{k+1}$. Then
$f_\t s_k=\sum_{\s\simk\t} s_{\t\s} (k) f_\s$ with
$$s_{\t\s}(k)=\frac{e_{\t\s}(k)-\delta_{\t\s}}{c_\t(k)+c_\s(k)}.$$
\end{enumerate}
\end{Theorem}


The following result gives  an explicit construction on primitive
idempotents and central primitive idempotents for $\cba{n}$. This
gives the  Wedderburn-Artin decomposition for  Brauer algebras
$\cba{n}$. Such results can be proved by the arguments in
\cite{M:gendeg}.  In \cite[3.16]{M:semi}, Mathas proved this result
for  a class of cellular algebras satisfying a ``separated
condition", which is analogue to the assumption~3.1. In our case,
such idempotents can be computed explicitly since, in the next
section,  we give recursive formulae for $\langle f_\t, f_\t\rangle$
for all $\t\in \UPD_n(\lambda)$ and all $(f, \lambda)\in \Lambda_n$.

\begin{Prop}\label{wed} Let $\cba{n}$ be a Brauer
algebra over a field $F$.
\begin{enumerate}\item Suppose $\t\in \UPD_n(\lambda)$. Then
$\frac{1}{\langle f_\t, f_\t \rangle} f_{\t\t}$ is a  primitive
idempotent  of  $\cba{n}$ with respect to the cell module $\Delta(f,
\lambda)$.
\item $\sum_{\t\in \UPD_n(\lambda)} \frac{1}{\langle f_\t, f_\t \rangle}
f_{\t\t}$ is a central primitive idempotent. Furthermore,
$$\sum_{(\frac{n-|\lambda|} 2, \lambda)\in \Lambda_n}
 \sum_{\t\in \UPD_n(\lambda)} \frac{1}{\langle f_\t,
f_\t \rangle} f_{\t\t}=1.$$
\end{enumerate}
\end{Prop}

\section {Discriminants of  Gram matrices }
In this section, we compute Gram  determinants  associated to all
cell modules of $\cba{n}$. Unless otherwise stated, we assume that
$F$ is a field which satisfies the assumption~\ref{assumption}.
When we use  Lemma~\ref{ek}(a) in 
Proposition~\ref{downlowest}--\ref{key3}, we have to  assume that \cite[4.12]{NW} for $\cba{n}$
holds.

Gram determinants for  all cell modules for Hecke algebras of type
$A$ were computed by Dipper and James in \cite{DJ}. When $q=1$, they
are the formulae for symmetric groups $\mathfrak S_n$. Therefore,
our formulae for $\det G_{0,\lambda}$ are special cases of those
formulae.

\begin{Defn} Given  $\s\in \UPD_n(\lambda)$ with $\s_{n-1}=\mu$, define
 \begin{enumerate} \item $\hat\s\in
\UPD_{n-1}(\mu)$,  such that $\hat \s_j=\s_j$, for all $1\le j\le
n-1$,  \item  $\tilde \s\in \UPD_n(\lambda)$ such that $\tilde
\s_i=\t^\mu_i$, $1\le i\le n-1$ and $\tilde\s_n=\lambda$.
\end{enumerate}\end{Defn}

For any partition $\lambda=(\lambda_1, \lambda_2, \cdots, )$, let
$\lambda!=\prod_i \lambda_i !\in \mathbb N$.
For any $\t\in \UPD_n(\lambda)$, define 
 $$F_{\t, k}=
\prod_{\substack{r\in \mathscr R(k)\\
r\neq c_{{\t}}(k)}}    \frac {x_k-r}{c_{\t}(k)-r},\quad  1\le k\le
n.$$


\begin{Prop}\label{inductionkey} Assume that  $\t\in \UPD_n(\lambda)$ with
$(f, \lambda)\in \Lambda_n$. If $\t_{n-1}=\mu$ with $(l, \mu)\in
\Lambda_{n-1}$,  then $ \langle f_\t,\ f_\t \rangle =\frac {1}
{\delta^{l} \mu!} \langle f_{\hat\t},\ f_{\hat\t} \rangle {\langle
f_{\tilde \t},\ f_{\tilde \t}\rangle} $.
\end{Prop}

\begin{proof} By the definition of $F_\t$, we have
$F_\t =F_{\hat\t} F_{\t, n}$.   Using the definition of $\m_\t$, we have
$$f_{\t^\lambda \t} = F_{\t^\lambda}h
\m_{\t^\mu,\hat \t} F_{\hat \t} F_{\t, n}
$$
 for some
$h\in \cba{n}$. Note that $$\m_{\t^\mu,\hat \t} F_{\hat \t} F_{\hat
\t} \m_{\hat{\t}, \t^\mu}\equiv\langle f_{\hat \t}, f_{\hat\t}
\rangle \m_{\t^\mu, \t^{\mu}} \pmod {\cba{n-1}^{\vartriangleright
(l, \mu)}}.$$ By Theorem~\ref{ud}(c) and Theorem~\ref{xproduct} and
Corollary~\ref{f-filtration},
$$\begin{aligned} f_{\t^\lambda\t} f_{\t\t^\lambda}
\equiv  & F_{\t^{\lambda}}   \m_{\t} F_{\hat \t} F_{\hat \t}
\m_{\hat \t, \t^\mu} F_{\t, n}^2 \sigma (h)
F_{\t^{\lambda}}  \pmod {\cba{n}^{\vartriangleright (f, \lambda)}} \\
\equiv & \langle f_{\hat\t}, f_{\hat\t} \rangle F_{\t^{\lambda}}
\m_{\t^\lambda, \tilde\t} F_{\t, n}^2
\sigma(h)  F_{\t^\lambda} \pmod {\cba{n}^{\vartriangleright (f, \lambda)}} \\
\equiv & \frac{ \langle f_{\hat\t}, f_{\hat\t} \rangle } {\delta^l
\mu!}  F_{\t^{\lambda}} h F_{\t, n}  \m_\mu F_{\t^{\mu}}
F_{\t^{\mu}} \m_\mu F_{\t, n}
\sigma(h)  F_{\t^\lambda} \pmod {\cba{n}^{\vartriangleright (f, \lambda)}} \\
\equiv &  \frac{ \langle f_{\hat\t}, f_{\hat\t} \rangle }
{\delta^l \mu!} f_{\t^\lambda \tilde\t } f_{\tilde\t \t^{\lambda}}
 \pmod {\cba{n}^{\vartriangleright (f, \lambda)}}\\
\equiv & \frac{\langle f_{\hat\t}, f_{\hat\t} \rangle \langle
f_{\tilde\t}, f_{\tilde\t} \rangle } {\delta^l \mu!} f_{\t^\lambda
\t^\lambda} \pmod {\cba{n}^{\vartriangleright (f, \lambda)}} .\\
\end{aligned}$$
Consequently, we have the formula for $
 \langle f_\t, f_\t  \rangle$, as required.
\end{proof}

\begin{Lemma} \label{fs} Assume $\t\in \UPD_n(\lambda)$ with  $(f, \lambda)\in \Lambda_n$.
 If  $\t s_{k}\in \UPD_n(\lambda)$ with $\t
s_{k}\vartriangleleft \t$,   then $\langle f_{\t s_{k}}, f_{\t
s_{k}} \rangle =(1-(c_\t(k+1)-c_\t(k))^{-2}) \langle f_\t, f_\t
\rangle$.
\end{Lemma}
\begin{proof}
The result follows from Lemma~\ref{fss} and  $\langle f_{\t},
f_{\t} \rangle =\langle f_{\t}s_{k} , f_{\t} s_{k}\rangle$.
\end{proof}

\def\R{\mathscr R(\lambda)^{<p}}
\def\A{\mathscr A(\lambda)^{<p}}
\def\RR{\mathscr R(\mu)^{<p}}
\def\AA{\mathscr A(\mu)^{<p}}
\def\RB{\mathscr {AR}(\lambda)^{\ge p}}
\def\RRB{\mathscr {AR}(\mu)^{\ge p}}
\def\AAB{\mathscr A(\mu)^{\ge p}}

\begin{Defn} For any $\lambda\vdash n-2f$, let $ \mathscr
A(\lambda)$ (resp. $\mathscr R(\lambda))$ be the set of all
addable (resp. removable) nodes of $\lambda$. Given a removable
(resp. an  addable) node  $p=(k, \lambda_k)$ (resp. $(k,
\lambda_k+1)$) of $\lambda$, define
\begin{enumerate}
\item  $\mathscr R(\lambda)^{<p} =\{(l, \lambda_l)\in \mathscr
R(\lambda) \mid l>k\}$,
\item $\mathscr A(\lambda)^{<p} =\{(l, \lambda_l+1)\in \mathscr
A(\lambda) \mid l>k\}$,
\item $\mathscr {AR}(\lambda)^{\ge p} =\{(l, \lambda_l)\in \mathscr
R(\lambda) \mid l\le k\}\cup \{(l, \lambda_l+1)\in \mathscr
A(\lambda) \mid l\le  k\}$.
\end{enumerate}
\end{Defn}

\begin{Prop}\label{key0}
Suppose  $\t\in \UPD_n(\lambda)$ with $\lambda\vdash n-2f$. If
$\hat \t=\t^\mu$ and   $\t_{n}= \t_{n-1}\cup \{p\}$ with $p=(k,
\lambda_k)$, then
\begin{equation}\label{key1}
{\langle f_\t,  f_\t \rangle}=- {\delta^f \mu!} \frac{ \prod_{q\in
\A} ( c_{\lambda} (p)+c_{\lambda} (q))}{\prod_{r\in \R} (c_{\lambda}
(p)-c_{\lambda} (r))}.
\end{equation}
\end{Prop}

\begin{proof} By assumption, $\t=\t^\lambda s_{a_k, n}$,  where
$a_k=2f+\sum_{i=1}^k \lambda_i$. Using Lemma~\ref{fs} repeatedly for
the pairs $\{f_{\t^\lambda s_{a_k, j}}, f_{\t^\lambda s_{a_k,
j+1}}\}$ with $a_k\le j\le n-1$, and noting that $\t \lhd \t s_{n-1}\cdots\lhd \t s_{n,
a_k}=\t^{\lambda}$,  we have $$\langle f_\t,  f_\t \rangle = \langle
f_{\t^{\lambda}}, f_{\t^{\lambda}}\rangle\prod_{j=a_k}^{n-1} (1-
(c_{\t} (n)-c_\t(j))^{-2}).$$

Since $f_{\t^\lambda}\equiv \m_{\lambda} \pmod
{\cba{n}^{\vartriangleright (f, \lambda)}}$, $\langle
f_{\t^{\lambda}}, f_{\t^{\lambda}}\rangle=\delta^f \lambda!$.   If
$a, b\in \mathbb Z$ with $b>a$, then
\begin{equation}\label{identity1}
\prod_{i=a}^b \frac { (q-i)^2 -1}{(q-i)^2}  =\frac{q-a+1} {q-a}
\frac{q-b-1}{q-b}\end{equation} We compute  $\prod_{j=a_k}^{n-1} (1-
(c_{\t} (n)-c_\t(j))^{-2})$ along each row of $\ts^{\lambda}$ first.
Via (\ref{identity1}), we need only consider the first and the last
nodes in each row. Note that $c_{\lambda}(p)=-c_\t(n)$. By
(\ref{identity1}) again, we have  $$ \prod_{j=a_k}^{n-1} (1- (c_{\t}
(n)-c_\t(j))^{-2}) =-\frac{1}{\lambda_k} \frac{ \prod_{q\in \A}
(c_{\lambda} (p)+c_{\lambda} (q))}{\prod_{r\in \R} (c_{\lambda}
(p)-c_{\lambda} (r))}, $$ proving (\ref{key1}).\end{proof}

\begin{Prop} \label{downlowest} Suppose $\t \in \UPD_n(\lambda)$ with
$\lambda=(\lambda_1, \dots, \lambda_k)\vdash n-2f$. If $\t^{\mu}
=\hat \t$ and $\t_{n-1} = \t_n \cup p$ with $p=(k,\mu_k)$, then
$${ \langle f_\t,\ f_\t\rangle }=A\cdot {\delta^{f-1}\mu!}\cdot (\delta+2\mu_k-2k)
\prod_{\substack{q\neq p\\ q\in \RB}}
\frac{c_\lambda(p)+c_\lambda(q)}{c_\lambda(p)-c_\lambda(q)}$$ where
$A=1$ if $\lambda_k=0$ and $A=\delta+\lambda_k-2k$ if $\lambda_k>0$.
 \end{Prop}

\begin{proof} Let $a = 2(f-1)+\sum_{j=1}^{k-1} \mu_j+ 1$.
Write $\tilde \m_\lambda=\m_{\t^\lambda}$ for $\t^\lambda\in
\UPD_{n-2}(\lambda)$. We have
$$\begin{aligned} f_\t e_{n-1} & \equiv
e_{2f-1}s_{2f,n}s_{2f-1,n-1}\tilde \m_{\lambda}\sum_{j=a}^{n-1}
s_{n-1,j}F_\t e_{n-1}  \pmod {\cba{n}^{\vartriangleright(f, \lambda)}} \\
& \equiv  e_{2f-1}s_{2f,n}s_{2f-1,n-1} \tilde
\m_{\lambda}\sum_{j=a}^{n-1}
s_{n-1,j}F_{\t,n-1}F_{\t,n}e_{n-1}\prod_{k=1}^{n-2} F_{\t, k}\pmod
{\cba{n}^{\vartriangleright(f,
\lambda)}}\\
\end{aligned}$$

Since $e_{n-1} x_{n-1}^k e_{n-1}=\delta_{n-1}^{(k)}$ for some
$\delta_{n-1}^{(k)}\in F[x_1, x_2, \cdots, x_{n-2}]\cap
Z(\cba{n-2})$,
$$e_{n-1}F_{\t,n-1}F_{\t, n}e_{n-1}=\Phi_\t(x_1,\cdots,
x_{n-2})e_{n-1},$$ for some  $\Phi_\t (x_1,\cdots, x_{n-2})\in
F[x_1, x_2, \cdots, x_{n-2}]\cap Z(\cba{n-2})$. In
\cite[2.4]{Nazarov:brauer}, Nazarov proved that
$$x_k^i s_k= s_k x_{k+1}^i  +\sum_{j=1}^{i} x_k^{i-j} (e_k-1)
x_{k+1}^{j-1} $$ for $\cba{n}$ over 
$\mathbb C$. In fact, it holds  for $\cba{n}$ over $F$ \cite[2.3]{NW}.    
 Therefore, using the above equality for $k=n-2$,
we have
$$e_{n-1}s_{n-2}F_{\t,n-1}F_{\t,n}e_{n-1}= \Psi_\t(x_1,\cdots,
x_{n-2})e_{n-1}
$$ for some  $ \Psi_\t(x_1,\cdots, x_{n-2})\in F[x_1,
x_2, \cdots, x_{n-2}]$. In order to simplify the notation, we use
$\Phi_\t$ and $\Psi_\t$ instead of $\Phi_\t(x_1,\cdots, x_{n-2})$
and $\Psi_\t(x_1,\cdots, x_{n-2})$ respectively. Let
$\Phi_{\t,\lambda}$ (resp. $\Psi_{\t, \lambda}$) be obtained  by
using $c_{\t^{\lambda}}(k)$ instead of $x_k$ in $\Phi_\t$
(resp.$\Psi_\t$). Then
$$\begin{aligned} & e_{2f-1} \cdots e_{n-2}\tilde \m_{\lambda}\Phi_\t
e_{n-1} \equiv \Phi_{\t,\lambda} \m_\u+\sum_{\substack {\v\in
\UPD_n(\lambda)\\ \v_{n-1}=\u_{n-1}\\
 \v_{n-2}\vartriangleright
\u_{n-2}}} b_\v\m_\v \pmod {\cba{n}^{\vartriangleright (f, \lambda)}}
\\
& e_{2f-1} \cdots e_{n-2}\tilde \m_{\lambda}\Psi_\t
e_{n-1}\sum_{j=a}^{n-2}
s_{n-2,j}  \equiv(\mu_k-1)\Psi_{\t,\lambda} \m_\u\\
& \hskip5.9cm +\sum_{\substack {\v\in
\UPD_n(\lambda)\\ \v_{n-1}=\u_{n-1}\\
 \v_{n-2}\vartriangleright
\u_{n-2}}}
a_\v \m_\v \pmod {\cba{n}^{\vartriangleright (f, \lambda)}}, \\
\end{aligned}
$$
where $\u\simnn \t$ with
$$\u_{n-1}=\begin{cases} \t_{n-2}\cup \{(k+1,1)\}, & \text{if
$\mu_k>1,$}\\ \t_{n-1}, &\text{if $\mu_k = 1$.}\end{cases}
$$
We remark that $\m_\u=e_{2f-1}\cdots e_{n-2} e_{n-1} \tilde
\m_\lambda$ and  $\u\vartriangleleft \v$ for any $\v\in
\UPD_n(\lambda)$ with $\v\simnn \t$. Therefore,
$$f_\t e_{n-1}\equiv (\Phi_{\t,\lambda}+(\mu_k-1)\Psi_{\t,\lambda}) \m_\u +\sum_{ \substack
{\v\in \UPD_n(\lambda)\\ \v_{n-1}=\u_{n-1}\\ \v_{n-2}
\vartriangleright \u_{n-2}}}(a_\v+b_\v)\m_\v \prod_{k=1}^{n-2}
F_{\t, k} \pmod {\cba{n}^{\vartriangleright (f, \lambda)}}.$$ In
particular, $e_{\t\u}(n-1)=\Phi_{\t,\lambda}+(\mu_k-1)\Psi_{\t,
\lambda}$. On the other hand, we have 
$$\begin{aligned}  \Phi_{\t, \lambda} f_\t e_{n-1}
\equiv & f_\t
e_{n-1}F_{\t,n-1}F_{\t,n}e_{n-1}\pmod {\cba{n}^{\vartriangleright (f, \lambda)}}\\
\equiv & \sum_{\v\simnn\t} e_{\t\v}(n-1) f_\v
F_{\t,n-1}F_{\t,n}e_{n-1}\pmod {\cba{n}^{\vartriangleright (f,
\lambda)}}\\
\equiv &
e_{\t\t}(n-1)f_\t e_{n-1}\pmod {\cba{n}^{\vartriangleright (f, \lambda)}}.\\
\end{aligned}
$$
Since we are assuming that \cite[4.12]{NW} holds, $e_{\t\t}(n-1)\neq 0$. So
$f_\t e_{n-1}\neq 0$ and $ e_{\t\t}(n-1)=\Phi_{\t, \lambda}$. If $\mu_k=1$, then
$(\mu_k-1)\Psi_{\t, \lambda}=0$. Suppose
 $\mu_k>1$. Then  $s_{n-2}\in \mathfrak S_\mu$. If we use $\Psi_{\t,
\lambda}$ instead of $\Phi_{\t, \lambda}$ above,  we obtain $
e_{\t\t}(n-1)=\Psi_{\t, \lambda}$. 
In this case, $f_\t s_{n-2}=f_\t$ since  both $\t_{n-1}\ominus \t_{n-2}$ and
$\t_{n-2}\ominus\t_{n-3}$  are in the same row.
Therefore, $$ \mu_k
e_{\t\t}(n-1)=\Phi_{\t, \lambda}+(\mu_k-1)\Psi_{\t,
\lambda}=e_{\t\u}(n-1).$$

Now, we compute $\langle f_\u, f_\u\rangle$.  By the similar
argument as above, we have
$$
\begin{aligned}
& f_{\t^{\lambda}\u}f_{\u \t^{\lambda}}\\ \equiv &
F_{\t^{\lambda}}e_{2f-1}s_{2f,n}s_{2f-1,n-1}\tilde  \m_{\lambda}
F_{\u,n-1}F_{\u,n}\tilde
\m_{\lambda}s_{n-1,2f-1}s_{n,2f}e_{2f-1}F_{\t^{\lambda}} \pmod
{\cba{n}^{\vartriangleright(f, \lambda)}}\\  \equiv &
F_{\t^{\lambda}} e_{2f-1} \cdots e_{n-2}\tilde  \m_{\lambda}
e_{n-1}F_{\u,n-1}F_{\u,n}e_{n-1}\tilde  \m_{\lambda}
e_{n-2} \cdots e_{2f-1}F_{\t^{\lambda}} \pmod {\cba{n}^{\vartriangleright(f, \lambda)}}\\
\equiv & F_{\t^{\lambda}} e_{2f-1} \cdots e_{n-2}\tilde \m_{\lambda}
e_{n-1}\Phi_\u \tilde  \m_{\lambda} e_{n-2} \cdots
e_{2f-1}F_{\t^{\lambda}} \pmod {\cba{n}^{\vartriangleright(f, \lambda)}}\\
\equiv & \delta^{f-1} \lambda !  F_{\t^{\lambda}} e_{2f-1} \cdots
e_{n-2} e_{n-1}\tilde  \m_{\lambda} \Phi_\u e_{n-2} \cdots
e_{2f-1}F_{\t^{\lambda}} \pmod {\cba{n}^{\vartriangleright(f,
\lambda)}}
\\
\end{aligned}
$$
Note that  $\Phi_\u\in F[x_1, \cdots, x_{n-2}]$. By
Lemma~\ref{fxproduct} and  Corollary~\ref{f-filtration},
$$\begin{aligned} & F_{\t^{\lambda}} e_{2f-1} \cdots e_{n-2} e_{n-1} \tilde \m_{\lambda}
\Phi_\u \\ \equiv & \Phi_{\u, \lambda} F_{\t^{\lambda}} e_{2f-1}
\cdots e_{n-2}
e_{n-1} \tilde \m_{\lambda} +\sum_{\substack{\v\in \UPD_n(\lambda)\\
\v_{n-1}=\u_{n-1}\\ \v_{n-2}\vartriangleright\u_{n-2}}} a_\v
F_{\t^\lambda} f_\v
\pmod {\cba{n}^{\vartriangleright (f, \lambda)}}\\
\end{aligned} $$
By Lemma~\ref{epro}, $f_\v e_{n-2} \cdots e_{2f-1}\pmod
{\cba{n}^{\vartriangleright (f, \lambda)}}$ can be written as a
linear combination of $f_\s$ with $\s_{n-1}=\v_{\n-1}$ which is not
equal to $(\t^{\lambda})_{n-1}$ under the assumption
$\t_{n-2}=\lambda\neq \varnothing$. If $\t_{n-2}=\varnothing$, then
there is no $\v\in \UPD_n(\lambda)$ such that
$\v_{n-2}\rhd\u_{n-2}$. In any case, for $\v\in \UPD_n(\lambda),$
$\v_{n-1}=\u_{n-1}$ and $\v_{n-2}\vartriangleright\u_{n-2}$, we
have, by Lemma~\ref{fst}, that
$$f_\v
e_{n-2} \cdots e_{2f-1} F_{\t^{\lambda}} \equiv 0\pmod
{\cba{n}^{\vartriangleright (f, \lambda)}}.$$ Therefore,
$$\begin{aligned} f_{\t^{\lambda}\u}f_{\u \t^{\lambda}}
& \equiv \Phi_{\u, \lambda} \delta^{f-1} \lambda! F_{\t^{\lambda}}
e_{2f-1}\cdots e_{n-1}\tilde  \m_{\lambda} e_{n-2}\cdots e_{2f-1}
F_{\t^{\lambda}}\pmod {\cba{n}^{\vartriangleright (f, \lambda)}}\\
& \equiv  \Phi_{\u, \lambda} \delta^{f-1} \lambda !
f_{\t^{\lambda}{\t^\lambda}} \pmod {\cba{n}^{\vartriangleright(f, \lambda)}}, \\
\end{aligned}
$$
forcing  $\langle f_\u, f_u\rangle = e_{\u\u}(n-1) \delta^{f-1}
\lambda !$. By Theorem~\ref{gram}(c)
$$ \frac{\langle f_\t, f_\t\rangle }{
\delta^{f-1}\mu!}=  \frac{ 1 }{\delta^{f-1}\mu!}\frac{
\mu_k^2e_{\t\t}(n-1)}{ e_{\u\u}(n-1)}e_{\u\u}(n-1) \delta^{f-1}
\lambda ! =\mu_k e_{\t\t}(n-1) .$$

Finally, we rewrite $e_{\t\t}(n-1)$ via Lemma~\ref{ek}a to obtain the
formulae, as required. 
\end{proof}

\begin{Prop}\label{key3}  Suppose $\t \in \UPD_n(\lambda)$ with
$(f, \lambda)\in \Lambda_n$, and $l(\lambda)=l$. If $\hat
\t=\t^{\mu}$, and $\t_{n-1} = \t_n \cup p$ with $p=(k,\mu_k)$ $k<
l$,  define $\u= \t s_{n,a+1}$ with  $a = 2(f-1)+\sum_{j=1}^{k}
\mu_j$ and $\v = (\u_1,\cdots, \u_{a+1})$. Define
$$
A= -\delta^{f-1} \mu! (\delta+2\mu_k-2k) \prod_{\substack{q\neq p\\
q\in\RB}}
\frac{c_\lambda(p)+c_\lambda(q)}{c_\lambda(p)-c_\lambda(q)}
\frac{\prod_{q\in \mathscr A(\mu)^{<p}} (c_{\mu}(p)-c_{\mu}(q))}
{\prod_{r\in \mathscr R(\mu)^{<p}} (c_{\mu}(p)+c_{\mu}(r))}$$

Then $$ {\langle f_\t,f_\t\rangle } =
\begin{cases} A, &
\text{if $(k, \lambda_k)\in R(\lambda)$}\\
\frac{A}{\delta-2+2\mu_k-2k}, &
\text{if $(k, \lambda_k)\not\in R(\lambda)$}\\
\end{cases}
$$
\end{Prop}

\begin{proof} By  Proposition~\ref{inductionkey}, Lemma~\ref{fs}, Proposition~\ref{downlowest} and
(\ref{key1}),  we have
$$\begin{aligned} &
{\langle f_\t,f_\t\rangle }= -\langle f_\u,f_\u \rangle
\frac{1}{\delta+\lambda_k-2k}\frac{\prod_{q\in \mathscr
A(\mu)^{<p}} (c_{\mu}(p)-c_{\mu}(q))}
{\prod_{r\in \mathscr R(\mu)^{<p}} (c_{\mu}(p)+c_{\mu}(r))}\\
& \langle f_\u,f_\u \rangle=\langle f_\v, f_\v \rangle
\lambda_{k+1}! \cdots \lambda_{l}!, \\
& \langle f_\v, f_\v \rangle=\delta^{f-1} \lambda_1 !\cdots
\lambda_k !\mu_k^2
e_{\v\v}(a).\\
\end{aligned}
$$ Combining Lemma~\ref{ek} and the above equalities, we obtain Proposition~\ref{key3}.
Note that  $\mu_k>1$ in this case.
\end{proof}

By  the classical branching rule for  $\cba{n}$ (see
\cite{Wenzl:ssbrauer}),
$$\Delta(f, \lambda)\downarrow \cong \bigoplus_{\mu\rightarrow
\lambda} \Delta(f, \mu)\bigoplus  \bigoplus_{\lambda\rightarrow
\nu} \Delta(f-1, \nu),$$ where $\Delta(f, \lambda)\downarrow$ is
the restriction of $\Delta(f, \lambda)$ to $\cba{n-1}$. We
 write $(l, \mu)\rightarrow (f, \lambda)$ if $\Delta(l, \mu)$
 appears in $\Delta(f,
\lambda)\downarrow $. Let $G_{l, \mu}$ be the Gram matrix associated
to the cell module $\Delta(l, \mu)$ which is defined by its
Jucys-Murphy basis. Let $\det G_{l, \mu}$  be the determinant of
$G_{l, \mu}$.

\begin{Defn} Suppose $(f, \lambda)\in \Lambda_n$ and $(l,
\mu)\in \Lambda_{n-1}$ such that $(l, \mu)\rightarrow (f, \lambda)$.
For any $\t\in \UPD_n(\lambda)$ with $\hat \t=\t^\mu\in
\UPD_{n-1}(\mu)$, define $\gamma_{\lambda/\mu}\in F$ to be the
scalar by declaring that
$$\gamma_{\lambda/\mu}= \frac{\langle f_\t,
f_\t\rangle}{\delta^l\mu!}. $$
\end{Defn}

Theorem~\ref{main} provides $\det G_{f, \lambda}$ only recursively.
F. L\"{u}beck wrote a GAP program for this recursive formula. It
takes us about two hours to compute $\det G_{f, \lambda}$ for all
$(f, \lambda)\in \Lambda_n$ with $n\le 35$  using a personal
computer.

\begin{Theorem}\label{main} Let $\cba{n}$ be a Brauer algebra over $\mathbb
Z[\delta]$. Let $\det G_{f, \lambda}$ be the Gram determinant
associated to the  cell module $\Delta(f, \lambda)$ of $\cba{n}$.
Then $$\det G_{f, \lambda}=\prod_{(l, \mu)\rightarrow (f,
\lambda)} \det G_{l, \mu} \cdot  \gamma_{\lambda/\mu}^{\dim
\Delta(l, \mu)}\in \mathbb Z[\delta].$$
 Furthermore, each scalar $\gamma_{\lambda/\mu}$ can be
computed explicitly  by Proposition~\ref{key0}, ~\ref{downlowest}
and ~\ref{key3}.
\end{Theorem}

\begin{proof} We first consider Brauer algebras over $\mathbb C(\delta)$,
where $\delta$ is an indeterminate. Obviously, $\mathbb C(\delta)$
satisfies the assumption~\ref{assumption}.
In order to  use the seminormal representation  constructed in \cite{NW}, 
we  consider $\cba{n}$ over the complex field. Note that $\mathbb R\subset \mathbb C$. 
By \cite[5.4a]{NW},   Lemma~\ref{ek}(a) holds  for 
$\cba{n}$ over $\mathbb C$ for infinite many $\delta$.
  Using the fundamental theorem 
of algebra, we obtain the result over $\mathbb C(\delta)$ where $\delta$ is an indeterminate. 
Therefore, we can use previous results in this section. Note that
the Gram matrix $\tilde G_{f, \lambda}$ which is defined via
orthogonal basis of $\Delta(f, \lambda)$ is a diagonal matrix. Each
diagonal is of the form $\langle f_\t, f_\t \rangle$, $\t\in
\UPD_n(\lambda)$. Therefore, $\det \tilde G_{f, \lambda}=\prod_{\t\in
\UPD_n(\lambda)} \langle f_\t, f_\t \rangle$. By Proposition~\ref{inductionkey} and
Corollary~\ref{equal},
$$\det G_{f, \lambda}=\det \tilde G_{f, \lambda}=\prod_{(l, \mu)\rightarrow (f,
\lambda)} \det G_{l, \mu} \cdot \gamma_{\lambda/\mu}^{\dim \Delta(l,
\mu)}.$$ Since the Jucys-Murphy basis of $\Delta(f, \lambda)$ is
defined over $\mathbb Z[\delta]$, the Gram matrices associated to
$\Delta(f, \lambda)$  which are  defined  over $\mathbb Z[\delta]$
and $\mathbb C(\delta)$ are same. We have  $\det G_{f, \lambda}\in
\mathbb Z[\delta]$ as required.
\end{proof}

The Gram determinant given in Theorem~\ref{main} is in $\mathbb
Z[\delta]$. Note that the Brauer algebra $\cba{n}_F$ over an
arbitrary field $F$ is isomorphic to $\cba{n}_{\mathbb
Z[\delta]}\otimes_{\mathbb Z[\delta]} F$. Therefore, one can get the
formula for the Gram determinant over an arbitrary field $F$ by
specialisation.

We can give a second proof of Theorem~\ref{semisimple} via
Theorem~\ref{main}.
 We will not give the details here. Instead, we will prove a result for
Birman-Murakami-Wenzl algebra in \cite{RS:ssBMW}, which is similar
to Theorem~\ref{main}. Via such a formula, we give a necessary and
sufficient condition for  Birman-Murakami-Wenzl algebras being
semisimple over an arbitrary field. Certain sufficient conditions
for the semisimplicity of Birman-Murakami-Wenzl algebras over
complex field  were obtained by Wenzl in \cite{W1}.

\begin {Example} The Gram determinant associated to   the cell module $\Delta(1,
\lambda)$ for  $\cba{4}$ with $\lambda=(2)\vdash 2$. \end{Example}

There are six elements in $\UPD_4(\lambda)$ as follows.
\unitlength 1mm 
\linethickness{0.4pt}
\ifx\plotpoint\undefined\newsavebox{\plotpoint}\fi 

\begin{picture}(87.87,70)(7,60)
\put(20,123){\circle{3}} \put(20.75,112.75){\circle{3}}
\put(21,102.75){\circle{3}} \put(21.5,92){\circle{3}}
\put(21.25,82){\circle{3}}
\put(29.25,120.75){\framebox(4.25,4)[cc]{}}
\put(29.5,110.75){\framebox(4.25,4)[cc]{}}
\put(29.5,101){\framebox(4,3.75)[cc]{}}
\put(29.75,89.75){\framebox(4.5,4)[cc]{}}
\put(30,80){\framebox(4.25,4)[cc]{}} \put(41.75,123.5){\circle{3}}
\put(49.25,121.25){\framebox(4,3.75)[cc]{}}
\put(59.25,121){\framebox(7,3.75)[cc]{}}
\put(49.5,111.25){\framebox(4.25,3.75)[cc]{}}
\put(38.5,111){\framebox(7.25,4)[cc]{}}
\put(59,111.25){\framebox(7.5,3.5)[cc]{}}
\put(40.5,100){\framebox(4.25,6.75)[cc]{}}
\put(50.25,100.75){\framebox(4,3.75)[cc]{}}
\put(59.5,100.75){\framebox(7.25,4)[cc]{}}
\put(38.5,90){\framebox(7.25,4)[cc]{}}
\put(48.5,90){\framebox(10.75,4)[cc]{}}
\put(38.5,80.25){\framebox(7.25,3.75)[cc]{}}
\put(48.5,80.5){\framebox(8.5,4)[cc]{}}
\put(60.5,80.5){\framebox(7.75,4)[cc]{}}
\put(21.5,72.5){\circle{3.61}}
\put(30.5,69.75){\framebox(4,4)[cc]{}}
\put(41.25,67.25){\framebox(4.25,8.25)[cc]{}}
\put(49,71){\framebox(7.75,4)[cc]{}}
\put(62.75,125){\line(0,-1){4}} \put(42.5,115){\line(0,-1){4}}
\put(63.25,114.75){\line(0,-1){4}} \put(40.75,103){\line(1,0){4}}
\put(63.5,104.5){\line(0,-1){3.5}}
\put(41.75,93.75){\line(0,-1){4.25}}
\put(52,93.75){\line(0,-1){3.75}} \put(52,90){\line(0,1){0}}
\put(55.75,93.75){\line(0,-1){3.75}}
\put(42,83.75){\line(0,-1){3.75}}
\put(48.5,80.25){\line(0,-1){3.25}} \put(53,84.25){\line(0,-1){7}}
\put(49,77.25){\line(1,0){4}} \put(64.75,84.5){\line(0,-1){4}}
\put(41.75,71.75){\line(1,0){3.75}}
\put(49,71.25){\line(0,-1){3.75}} \put(49,67.5){\line(1,0){4.5}}
\put(61,90){\framebox(7,3.75)[cc]{}}
\put(60.75,69.75){\framebox(7.5,4.25)[cc]{}}
\put(64.5,93.5){\line(0,-1){3.5}} \put(64.5,74){\line(0,-1){4}}
\multiput(22.5,126)(-.03373016,-.04365079){126}{\line(0,-1){.04365079}}
\multiput(43.5,126.25)(-.03365385,-.05528846){104}{\line(0,-1){.05528846}}
\multiput(22.75,115.25)(-.03361345,-.03991597){119}{\line(0,-1){.03991597}}
\multiput(23.25,105.5)(-.03373016,-.03968254){126}{\line(0,-1){.03968254}}
\multiput(23.75,94.75)(-.03368794,-.03900709){141}{\line(0,-1){.03900709}}
\multiput(23.5,84.75)(-.03373016,-.04166667){126}{\line(0,-1){.04166667}}
\multiput(24,75.5)(-.03368794,-.04255319){141}{\line(0,-1){.04255319}}
\put(53,74.75){\line(0,-1){7.25}} \put(7.5,122.5){$\t_1$}
\put(7.5,112.25){$\t_2$} \put(7.5,102.75){$\t_3$}
\put(7.5,91.75){$\t_4$} \put(7.5,80.75){$\t_5$}
\put(7.5,71.5){$\t_6$}\put(15,112.25){(} \put(15,122.5){( }
\put(15,102.75){ ( } \put(15,91.75){ ( } \put(15,80.75){ ( }
\put(15,71.5){ ( } \put(70,112.25){ ) } \put(70,122.5){ ) }
\put(70,102.75){ ) } \put(70,91.75){ ) } \put(70,80.75){ ) }
\put(70,71.5){ ) }

\put(87,123.4){$\langle f_{\t_1}, f_{\t_1}\rangle = \langle
f_{\hat{\t}_1}, f_{\hat{\t}_1}\rangle 2$} \put(87,112.89){$\langle
f_{\t_2}, f_{\t_2}\rangle = \langle f_{\hat{\t}_2},
f_{\hat{\t}_2}\rangle  2$} \put(87,102.17){$\langle f_{\t_3},
f_{\t_3}\rangle =\langle f_{\hat{\t}_3}, f_{\hat{\t}_3}\rangle  2$}
\put(87,91.45){$ \langle f_{\t_4}, f_{\t_4}\rangle = \langle
f_{\hat{\t}_4}, f_{\hat{\t}_4}\rangle  3e_{\t_4\t_4}(3)$}
\put(87,82.41){$\langle f_{\t_5}, f_{t_5}\rangle = \langle
f_{\hat{\t}_5}, f_{\hat{\t}_5}\rangle e_{\t_5\t_5}(3)$}
\put(87,71.27){$\langle f_{\t_6}, f_{\t_6}\rangle = \langle
f_{\hat{\t}_6}, f_{\hat{\t}_6}\rangle  e_{\t_5\t_5}(3)$}
\end{picture}

Define
 $\m_{\t_1} = e_1(1+s_3)$, $\m_{\t_2} =
e_1(1+s_3)s_2(1+s_1)$, $\m_{\t_3} = e_1(1+s_3)s_2s_1$, $ \m_{\t_4}
= e_1(1+s_3)s_2s_3s_1s_2(1+s_2+s_2s_1)$, $\m_{\t_5} =
e_1(1+s_3)s_2s_3s_1s_2$ and $\m_{\t_6} = e_1(1+s_3)s_2s_3s_1$.
Then $\{\m_{\t_i}\pmod  {\cba{4}^{\vartriangleright \lambda}} \mid
1\le i\le 6\}$ is the Jucys-Murphy base of $\Delta(1, (2))$. The
corresponding Gram matrix $G_{1,(2)}$  is given as follows:
$$
G_{1,(2)}=\left(\begin{array}{cccccc}
 2 \delta  & 4 & 2 & 4 & 0 & 2 \\
4 & 4 \delta + 4 & 2\delta + 2 & 8 & 4 & 2\\
2 & 2\delta + 2 & 2\delta & 4 & 2 & 2 \\
4 & 8 & 4 & 6\delta + 12 & 2\delta + 4 & 2\delta + 4 \\
0 & 4 & 2 & 2\delta + 4 & 2 \delta & 2 \\
2 & 2 & 2 & 2\delta + 4 & 2 & 2\delta
\end{array}\right)
 $$
Therefore, $\det
 G_{1,(2)} = 64{\delta}^3{(\delta - 2)}^2(\delta +
 4)$.
On the other hand, by the formulae on $\gamma_{\lambda/\mu}$,
 we obtain  $e_{\t_4\t_4}(3) =
\frac{\delta (\delta + 4)}{3(\delta + 2)}$ and $e_{\t_5\t_5}(3) =
\frac{2\delta (\delta - 2)}{3(\delta - 1)}$. By Theorem~\ref{main}
again, we recover  the formula for $\det
 G_{1,(2)}$.

 \providecommand{\bysame}{\leavevmode ---\ }
\providecommand{\og}{``} \providecommand{\fg}{''}
\providecommand{\smfandname}{and}
\providecommand{\smfedsname}{\'eds.}
\providecommand{\smfedname}{\'ed.}
\providecommand{\smfmastersthesisname}{M\'emoire}
\providecommand{\smfphdthesisname}{Th\`ese}

\end{document}